\documentclass{amsart}
\usepackage{graphicx}
\usepackage{color}
\usepackage{multirow}
\usepackage{array}
\usepackage{tikz}
\usetikzlibrary{spy}
\usetikzlibrary{backgrounds}
\usetikzlibrary{decorations}
\usetikzlibrary{shapes,snakes}
\usetikzlibrary{shapes.geometric}
\usetikzlibrary{shadows}
\tikzset{
   invisible/.style={opacity=0,text opacity=0},
   visible on/.style={alt={#1{}{invisible}}},
   alt/.code args={<#1>#2#3}{%
     \alt<#1>{\pgfkeysalso{#2}}{\pgfkeysalso{#3}} },}
\tikzset{cross/.style={cross out, draw=black, minimum size=2*(#1-\pgflinewidth), inner sep=0pt, outer sep=0pt},cross/.default={1pt}}
\def\d{\displaystyle}

\def\E{\mathbf{E}}
\def\P{\mathbf{P}}

\def\M{\mathcal{M}}

\def\barsmh{\bar{s}_{\hat{m}}}
\def\barsm{\bar{s}_m}
\def\barsmp{\bar{s}_{m'}}
\def\barga{\bar{\gamma}}
\def\hatsmh{\hat{s}_{\hat{m}}}
\def\hatsm{\hat{s}_m}
\def\hatsmp{\hat{s}_{m'}}
\def\mh{\hat{m}}

\def\Tb{\mathbf{T}}

\def\Eb{\mathbb{E}}

\def\Omf{\Omega_{m_f}(\varepsilon)}

\def\OO1{\Omega_1(\xi)}
\def\mO2{\Omega_2(\xi)}

\def\[[{[\![}
\def\O{\mathbf{1}_{\Omega(\varepsilon,\xi)}}
\def\Of1{\mathbf{1}_{\Omega_{m_f}(\varepsilon)}}

\newcommand{\Nu}{{\boldsymbol{\nu}}}
\newcommand{\T}{{\mathbf T}}
\newcommand{\barT}{{\bar{\mathbf T}}}
\newcommand{\barE}{{\bar{\mathbf E}}}

\newcommand{\thetabf}{\text{\mathversion{bold}{$\theta$}}}

\newtheorem{thrm}{Theorem}[section]
\newtheorem{lmm}[thrm]{Lemma}
\newtheorem{crllr}[thrm]{Corollary}
\newtheorem{prpstn}[thrm]{Proposition}

\newtheorem{rmrk}[thrm]{Remark}

\newcolumntype{M}[1]{>{\centering\arraybackslash}m{#1}}
\newcolumntype{N}{@{}m{0pt}@{}}

\begin{document}
\title{Model Selection for the Segmentation of Multiparameter Exponential Family Distributions}
\author{Alice Cleynen}\address{Harvard School of Public Health, Boston, USA}
\email{acleynen@jimmy.harvard.edu}
\author{Emilie Lebarbier}\address{AgroParisTech UMR518, Paris 5e, France}
\address{INRA UMR518, Paris 5e, France.}
\email{emilie.lebarbier@agroparistech.fr}


\begin{abstract}
We consider the segmentation problem of univariate distributions
from the exponential family with multiple parameters. In
segmentation, the choice of the number of segments remains a
difficult issue due to the discrete nature of the change-points. In
this general exponential family distribution framework, we propose a
penalized $\log$-likelihood estimator where the penalty is inspired
by papers of L. Birg\'e and P. Massart. The resulting estimator is
proved to satisfy an oracle inequality. We then further study the
particular case of categorical variables by comparing the values of
the key constants when derived from the specification of our general
approach and when obtained by working directly with the
characteristics of this distribution. Finally, a simulation study is
conducted to assess the performance of our criterion for the
exponential distribution,  and an application on real data modeled
by the categorical distribution is provided.
\end{abstract}

\maketitle

\section{Introduction}

Penalized-likelihood approaches are becoming more and more popular
in numerous domains in statistics: density estimation, variable
selection, machine learning, etc. Here we consider a multiple
change-point detection setting for univariate datasets where we are
interested in the estimation of the number of segments $K$. While
earlier works have focused on variables distributed from specified
distributions, for instance Gaussian
\cite{lebarbier_detecting_2005}, or either Poisson or negative
binomial distributions \cite{CleynenLebarbier2014}, we consider here
a more general framework of exponential family distributions. More
precisely, the model is the following:
\begin{eqnarray*}
Y_t \sim \mathcal{G}({\thetabf}_t) = s(t),\ 1\leq t\leq n
\end{eqnarray*}
where $\mathcal{G}$ is a distribution from the exponential family
with parameter ${\thetabf}_t \in \mathbb{R}^d$. In the context of
change-point models, we want to consider that $\thetabf_t$ is
piece-wise constant along the time-line $1,\dots,n$ and we therefore
wish to identify a partition of $\{1,\dots,n\}$ into $K$ segments
within which the observations can be modeled as following the same
distribution  while they differ between segments.

In our framework we will consider the minimal canonical form of the
exponential family distribution which we will write as
\begin{eqnarray*}
s(t)=h(Y_t)\exp\left[{\thetabf}_t.\T_t-A({\thetabf}_t) \right]
\end{eqnarray*}
where we use the $.$ symbol to denote the canonical scalar product of $\mathbb{R}^d$, $A$ is the log-partition function of $\mathcal{G}$ and $\T_t = T(Y_t)\ (\in \mathbb{R}^d)$ is the minimal sufficient statistic associated with variable $Y_t$. \\

While almost all methods for choosing the number of segments can be seen as penalized-likelihood approaches (Akaike Information Criterion, \cite{AIC}, Bayes Information Criterion \cite{Yao88}, Integrated Completed Likelihood \cite{rigaill2012exact}, etc), we and other authors (see for instance \cite{lebarbier_detecting_2005,BMMinPen,zhang2007modified}) have previously emphasized how crucial the choice of the penalty function is in contexts such as segmentation where the size of the collection of models grows with the size of the data. For this reason we have extended the approach developed in the  Poisson and negative binomial cases to more general distributions from the exponential family. The motivation for this work initially relied on the very strong similarities which we had observed between the Poisson and negative binomial distributions. Most of those similarities are in fact related to properties of the log-partition function $A$. \\

Exponential family distributions, and in particular the log-partition function, have been well studied in the past years. In a pioneer work \cite{brown1986fundamentals}, Brown has described the fundamental properties of exponential family distributions, such as parametrization using sufficient statistics, differentiability of the log-partition function and its relation to moments, etc. More recently, \cite{wainwright2008graphical} has demonstrated the strong links between graphical models and exponential family, \cite{kakade2009learning} has studied the sub-exponential growth of the cumulants of an exponential family distribution and studied convergence rates of regularization algorithm under sparsity assumptions while \cite{lee2013model} has studied consistency properties of the lasso procedure under some convexity assumption. \\

In this paper our goal is to mimic the procedure followed in the
Poisson and negative binomial cases for obtaining the oracle
inequality  while emphasizing the role of the log-partition
function. Considering the essential role played by the chi-square
statistic in this earlier work, we  will restrict the considered
families to those with positive marginal sufficient statistics
(hence allowing the definition of $\chi^2$), thus we will assume
that the set of natural parameters $\thetabf_t$ is restricted to a
compact $\boldsymbol{\nu}$ such that $\nabla A(\boldsymbol{\nu})$
belongs to $\{\mathbb{R}^+\}^d$. This for instance excludes the
Gaussian distribution unless we assume the mean parameter to be
positive.

Among the key features of minimal exponential family is the
relationship between the derivatives of $A$ and the moments of the
sufficient statistics. Hence the first two moments are given by
\begin{itemize}
\item $\mathbb{E}(\T_t)=\nabla A(\thetabf_t)$
\item $Cov(\T_t)= \nabla^2 A(\thetabf_t) $
\end{itemize}
Moreover, using minimal representation of the exponential family
ensures that the gradient mapping $\nabla A : \boldsymbol{\nu}
\rightarrow \nabla A(\Nu)$ is a bijection (see for instance
\cite{wainwright2008graphical}). These properties will be
fundamental in the construction of the oracle inequality. The proof
of the later mimicks the approach adopted in
\cite{CleynenLebarbier2014}, most of the time studying each marginal
of the sufficient statistic separately. \newline

In Section \ref{proc} we introduce the collection of models and the
penalized-likelihood framework and state our main result while in
Section \ref{exp} we derive exponential bounds for the sufficient
statistic. Proofs of our main statement and oracle inequality are
given in Section \ref{proof}. In Section \ref{classic} we explicit
the constants and assumptions used in a list of classical law from
the exponential family such as Poisson, exponential, gamma, beta,
etc, and we study the particular case of categorical variables in
Section \ref{categorical} to assess the  precision lost in dealing
with general exponential family instead of directly bounding the
particular distribution. Finally, in Section \ref{appli}, we
illustrate the performance of our approach on a simulation study
based on the exponential distribution, and propose an application to DNA
sequence distribution on a real data-set which we model using a piece-wise
constant categorical distribution.

\section{Model Selection Procedure} \label{proc}
\subsection{Penalized maximum-likelihood estimator}

In our change-point setting, we will want to consider partitions $m$
of the set $\{1,\dots,n\}$ on which our models will be piece-wise
constant. More precisely, for a given partition $m$ and denoting
$J$ a generic segment of $m$, we define the collection of models
associated to  $m$ as:
\begin{eqnarray*}
\mathcal{S}_m=\{s_m|\forall J \in m, \forall t \in J,
s_m(t)=\mathcal{G}(\thetabf_J) \}.
\end{eqnarray*}

We will consider the log-likelihood contrast $\gamma$ and the
associated  Kullback-Leibler distance
$K(s,u)=\mathbb{E}[\gamma(u)-\gamma(s)]$ between distribution $s$
and $u$ so that for distributions
 $s(t)=\mathcal{G}(\thetabf_t)=h(Y_t)\exp\left[{\thetabf}_t.\T_t-A({\thetabf}_t) \right]$, and $u=\left(\mathcal{G}(\boldsymbol{p}_t)\right)_t$, we have
\begin{eqnarray*}
\gamma(s)& = & \sum_{t=1}^n \left(A(\thetabf_t)-\thetabf_t.\T_t \right), \quad \text{and} \\
K(s,u) & = & \sum_{t=1}^n \left[ \nabla
A(\thetabf_t).(\thetabf_t-\boldsymbol{p}_t)-\left(A(\thetabf_t)-A(\boldsymbol{p}_t)
\right) \right].
\end{eqnarray*}

The minimal contrast estimator $\hatsm$ of $s$ on the collection
$\mathcal{S}_m$ is therefore
$\hatsm=\arg\min_{u\in\mathcal{S}_m}\gamma(u) $ and is given by
\begin{eqnarray*}
\hatsm(t)&=&h(Y_t)\exp\left[\nabla A^{-1}(\bar{\Tb}_J) .\Tb_t -A\left[\nabla A^{-1}(\bar{\Tb}_J)\right]\right],\\
\end{eqnarray*}
where $\Tb_J=\sum_{t\in J}\Tb_t$, is the sum of the sufficient statistics on segment $J$, and $\bar{\Tb}_J=\Tb_J/|J|$ its mean; the bijective mapping of the gradient of $A$ ensuring the existence and uniqueness of $\nabla A^{-1}(\bar{\T}_J)$. \\
Similarly, the projection $\barsm$ of $s$ in terms of
Kullback-Leibler on collection $\mathcal{S}_m$ is
$\barsm=\arg\min_{u\in\mathcal{S}_m} K(s,u) $ and is given by
\begin{eqnarray*}
\barsm(t)&=&h(Y_t)\exp\left[\nabla A^{-1}(\bar{\E}_J) .\Tb_t
-A\left[\nabla A^{-1}(\bar{\E}_J)\right]\right],
\end{eqnarray*}
where $\E_t=\nabla A(\thetabf_t)$, $\E_J=\sum_{t\in J}\E_t$ and
$\bar{\E}_J=\E_J/|J|$. \newline

As is classical in penalized-likelihood settings, since minimizing
the Kullback-Leibler risk would require knowing the true
distribution $s$,  we will wish to choose a penalty function
$pen(m)$ such that the penalized estimator $\hatsmh$, where $\hat{m}
= \arg\min \gamma(\hatsm)+pen(m)$, satisfies an oracle inequality of
the type $$\mathbb{E}[K(s,\hatsmh)]\leq
C_1\mathbb{E}[K(s,\hat{s}_{m(s)})] +C_2 $$ where $C_2$ is negligible
compared to $C_1\mathbb{E}[K(s,\hat{s}_{m(s)})]$.

To this effect, here as in previous works (see for example\cite{massart2007concentration}), we will introduce an event of large probability, $\Omf$ defined on a minimal partition $m_f$, where the fluctuation of each centered marginal is bounded. On this event we will derive tight controls of the risk of the models which will lead to the first part of the oracle inequality, \textit{i.e.} $\mathbb{E}[K(s,\hatsmh)]\leq C_1\mathbb{E}[K(s,\hat{s}_{m(s)})]$. On the complementary of this event, we will obtain less tight results which are compensated on expectation by the negligible  probability of the event. This control will result in the $C_2$ constant. \\
The choice of $\varepsilon$ is therefore crucial in insuring that both $C_1$ and $C_2$ are as small as possible, while having the negligibility property of $C_2$ compared to $C_1\mathbb{E}[K(s,\hat{s}_{m(s)})]$. In practice, this choice is data-driven and is efficiently performed through the use of the slope heuristic \cite{arlot2009data}. In this paper, we therefore consider a generic but fixed $\varepsilon$, and aim at obtaining the shape of the penalty function. This choice is given in the next section. \\

\subsection{Main result}

\begin{thrm} \label{theo}
Let $s=\{s(t)\}_{1\leq t\leq n}$ be a distribution from the
exponential family such that its parameters $\thetabf_t$ belong to a
compact set $\boldsymbol{\nu}$  with greatest band $\nu$ and such
that $\nabla A(\boldsymbol{\nu}) \subset \{\mathbb{R}^+\}^d$. Assume
there exist some positive constants $\zeta^i$ and $\kappa^i$ and a
positive constant $\kappa$ such that

\begin{eqnarray*}
 \forall 1\leq t \leq n,\ \forall 1\leq i \leq d,\quad \zeta^i E^i_t \leq Var T^i_t \leq \kappa^i E^i_t \quad \quad (\mathbf{H}) \\
 \text{and} \quad \kappa \geq \max_t \max_i \{\kappa^i, 1/R^i_t, \alpha_t^i\sqrt{Var T^i_t} \}
\end{eqnarray*}
where $R_t^i$ is the radius of convergence of $A(\thetabf_t+z\boldsymbol{\delta}_i)$, $\boldsymbol{\delta}_i=(0,\dots,0,1,0,\dots,0)$ and $\alpha_t^i$ is an exponential control of the growth rate of the coefficients of its power series. \\
Let $\zeta=min_i \{ \zeta_i\}$.  Let $\mathcal{M}_n$ be a collection
of partitions constructed on a partition $m_f$ such that there
exists $ \Gamma>0 $ satisfying $\forall J \in m_f, |J|\geq \Gamma
\log(n)^2$, and let $(L_m)_{m\in \mathcal{M}_n}$ be some family of
positive weights satisfying
\begin{eqnarray*}
     \Sigma = \sum_{m\in\mathcal{M}_n} \exp(-L_m|m|) < +\infty . \label{weights}
\end{eqnarray*}
Let $\varepsilon > 0$ and let $\beta_\varepsilon$ be a positive
constant depending on $\varepsilon$, distribution $\mathcal{G}$ and the observation $y_t$. If
for every $m \in \mathcal{M}_n$
   \begin{eqnarray*}
     pen(m)\geq \beta_{\varepsilon} d |m| \left(1 + 4\sqrt{L_{m}} \right)^2,
   \end{eqnarray*}
then
\begin{eqnarray} \label{Risk}
     \E\left[ h^2(s,\hatsmh) \right] &\leq & C_{\beta_{\varepsilon}} \inf_{m\in\mathcal{M}_n} \left\{K(s,\barsm) +pen(m)\right\} \nonumber \\
     &&+ C(d,\nu,\kappa,\Gamma,\beta_{\varepsilon},\Sigma,E^{max},E^{min}),
   \end{eqnarray}
where $E^{min} =\min_{1\leq i\leq d} \{\mathbb{E} (T^i) \} $ and $E^{max}= \max_{1\leq i\leq d}\{\mathbb{E} (T^i) \} $ are bounds on the expectation of the sufficient statistics.
\end{thrm}

\vspace{0.2cm}
\begin{rmrk} It can be noted that some assumptions of theorem \ref{theo} are in fact redundant and were included to simplify notations and clarify the dependency to the required assumptions. Indeed,
 \begin{itemize}
 \item the existence of $E^{min}$ and $E^{max}$ are guaranteed since we assume that $\thetabf$ belongs to a compact $\boldsymbol{\nu}$. Similarly, the greatest bound $\nu$ of the compact needs not be introduced and we could use $\theta^{max}-\theta^{min}$ instead. However in general the marginals of the sufficient statistic might not fluctuate in the same subset of $\mathbb{R}^+$ and therefore $\nu$ is a refinement of the bound.
 \item bounding $\thetabf$ (on $\Nu$) also implies that the variance of the sufficient statistics are bounded (since $A$ is analytic), and therefore the existence of  $\zeta$ and $\kappa$ in assumption $(\mathbf{H})$ are guaranteed.
 \end{itemize}
The theorem could therefore be re-written with the minimal
conditions that $\thetabf~\in~\Nu$, $\nabla A(\boldsymbol{\nu}) \subset \{\mathbb{R}^+\}^d$ and the hypothesis on $m_f$ and $\Sigma$.
\end{rmrk}

In our change-point setting, we will choose weights $L_m$ that
depend on $m$ only through their number of segments $|m|$, resulting
in $L_{m}=1.1+\log\left(\frac{n}{|m|} \right)$ (see
\cite{lebarbier_detecting_2005} for a justification of this choice).
This leads to a penalty function of the form:
\begin{equation}
pen(m)= \beta_{\varepsilon} d |m|
\left(1+4\sqrt{1.1+\log\left(\frac{n}{|m|}\right)} \right),
\label{penalty_exhaustive}
\end{equation}
and the inequality of  theorem \ref{theo} therefore becomes:
\begin{eqnarray*}
 \E\left[ h^2(s,\hatsmh) \right] &\leq & C_{\beta_{\varepsilon}} \inf_{m\in\mathcal{M}_n} \left\{K(s,\barsm) +d\beta_{\varepsilon}|m|\left(1+4\sqrt{1.1+\log\left(\frac{n}{|m|}\right)} \right)\right\} \\
&+&
C(d,\nu,\kappa,\Gamma,\beta_{\varepsilon},\Sigma,E^{max},E^{min}).
\end{eqnarray*}

\vspace{0.2cm} The following proposition gives a bound on the
Kullback-Leibler risk associated to $\hatsm$, which we prove in
Section \ref{risk}. This bound is then integrated into the main
theorem to obtain our oracle inequality in corollary
\ref{corollary}.
\begin{prpstn} \label{Risk-control}
Under the assumptions of theorem \ref{theo}, we have
\begin{eqnarray} \label{eq:Risk-control}
K(s,\barsm) + C(\varepsilon) d |m|-\dfrac{C(a,\Gamma,E^{min},E^{max},\varepsilon,\kappa)}{n^{a/2-1}}\leq\mathbb{E}[K(s,\hatsm)],
\end{eqnarray}
where $a\geq 2$ and  $C(\varepsilon)$ is a constant that depends on
$\varepsilon$, the constants involved in hypothesis $(\mathbf{H})$,
and the constraints on the compact $\Nu$.
\end{prpstn}
\vspace{0.2cm}
\begin{crllr} \label{corollary}
Let $s=\{s(t)\}_{1\leq t\leq n}$ be a distribution from the
exponential family such that its parameters $\thetabf_t$ belong to a
compact $\boldsymbol{\nu}$ of $\{\mathbb{R}^+\}^d$. Let
$\mathcal{M}_n$ be a collection of partitions constructed on a
partition $m_f$ such that there exists $ \Gamma>0 $ satisfying
$\forall J \in m_f, |J|\geq \Gamma \log(n)^2$, and assume $\exists \
\zeta$ and $\kappa$ such that assumption $(\mathbf{H})$ is verified.
There exists some constant $C$ such that
\begin{eqnarray} \label{oracle}
     \E\left[ h^2(s,\hatsmh) \right] &\leq & C \log(n) \inf_{m\in \M_n} \left\{\E[ K(s,\hatsm)]\right\} \nonumber \\
     &&+ C(d,\kappa,\nu,\Gamma,E^{min},E^{max},\beta_{\varepsilon},\Sigma).
   \end{eqnarray}

\end{crllr}

\section{Exponential bounds}\label{exp}
Following previous works
\cite{BMGaussian,castellan1999modifiedb,CleynenLebarbier2014}, we
will need to obtain exponential bounds on the fluctuation of the
variables in order to derive bounds on the risk of models. Starting
from equation
\begin{eqnarray*}
     \gamma(\hatsmh) + pen(\hat{m})& \leq& \gamma(\hatsm) + pen(m) \leq \gamma(\barsm) + pen(m),
   \end{eqnarray*}
and introducing the centered loss
$\barga(u)=\gamma(u)-\mathbb{E}[\gamma(u)]$, we write
\begin{eqnarray*}
      K(s,\hatsmh)& \leq&  K(s,\barsm) + \barga(\barsm) -\barga(\hatsmh) -pen(\mh) +pen(m).
   \end{eqnarray*}

We subsequently decompose $\barga(\barsm) -\barga(\hat{s}_{m'})$ for
any $m' \in \mathcal{M}_n$ in
   \begin{eqnarray} \label{decomposition}
     \barga(\barsm) -\barga(\hat{s}_{m'}) = \left(\barga(\bar{s}_{m'}) -\barga(\hat{s}_{m'})\right) + \left(\barga(s) - \barga(\bar{s}_{m'})\right) + \left(\barga(\barsm)
     -\barga(s)\right),
   \end{eqnarray}
to control each term separately, with the main term now defined on
the same model. This will lead us to introduce the chi-square
statistic
\begin{eqnarray}    \label{def-chi}
\chi^2_m \ =\  \sum_{i=1}^d \chi^2(m,i) \ = \ \sum_{i=1}^d
\sum_{J\in m}\dfrac{(T^i_J-E^i_J)^2}{E^i_J},
\end{eqnarray}
and control the fluctuations of the sufficient statistics from their
expectation.

\subsection{Control of $T_J$} \label{ControlYJ}

We are now interested in controlling the variables $T_J^i$ so that
we can later apply Bernstein's inequality to $\chi^2(m,i)$. To this
purpose, we wish to apply the large deviation result from Barraud
and Birg\'e (see lemma 3 in \cite{BaraudBirge}) which requires the
control of the Laplace transform of $T_J^i$.

Let $i\in \{1,\dots,d\}$ and $\boldsymbol{\delta}_i$ be the vector
of $\mathbb{R}^d$ with $1$ as the $i$th component and $0$ as others.
Let $1\leq t \leq n$ and let $z$ belong to $[-R_t^i,R_t^i]$, where
$R_t^i$ is the radius of convergence of the power series of
$A(\thetabf_t+z \boldsymbol{\delta}_i)$. We have
   \begin{eqnarray*}
     \log \Eb\left[e^{z(T_t^i-E_t^i)} \right] &=& \sum_{k\geq 2} c_k \dfrac{z^k}{k!} \ =\  c_2 \frac{z^2}{2} \sum_{k\geq 0} 2 \frac{c_{k+2}}{c_2} \dfrac{z^k}{(k+2)!},\\
   \end{eqnarray*}
where $c_k$ is the $k$th cumulant of random variable $T_t^i$. By
analyticity of $A$ on the natural parameter space (see Lemma 3.3
\cite{kakade2009learning} and \cite{brown1986fundamentals}), there
exists $\alpha_t^i>0$ $\forall 1\leq i\leq d$ such that
\begin{eqnarray*}
\forall k\geq 3, \ \left|\dfrac{c_k(T^i_t)}{c_2(T^i_t)^{k/2}}\right|
\leq \frac{1}{2} k! (\alpha_t^{i})^{k-2},
\end{eqnarray*}
and therefore if $\tau_i=\max_t \left\{\max\{ \alpha_t^i
\sqrt{Var(T^i_t)},1/R_t^i\} \right\}$, the Laplace transform of
$T_t^i$ can be bounded as
   \begin{eqnarray*}
     \log \Eb\left[e^{z(T_t^i-E_t^i)} \right] \leq    Var(T^i_t) \frac{z^2}{2} \sum_{k\geq 0}  (\tau_i z)^k \leq  \left\{\begin{array}{ll} \dfrac{z^2}{2}\dfrac{\kappa^i {E_t^i}}{1-\tau_i z} & \text{if } 0\leq z<1/\tau_i \\
     \dfrac{z^2}{2}\kappa^i{E_t^i} &  \text{if } -1/\tau_i<z<0 \end{array}\right. .
   \end{eqnarray*}
We can therefore apply the large deviation result from
\cite{BaraudBirge}:
   \begin{eqnarray*}
     P\left[\sum_{t\in J}  (T^i_t-E^i_t) \geq \sqrt{2\kappa^i x E^i_J}+\tau_i x \right] \leq e^{-x}\\
     P\left[-\sum_{t\in J}  (T^i_t-E^i_t) \leq \sqrt{2\kappa^i x E^i_J}\right] \leq e^{-x}
   \end{eqnarray*}
which implies, denoting $\kappa=\max_i\max\{\kappa^i,\tau_i \}$,
that
   \begin{eqnarray*}
     P\left[\left|(T^i_J-E^i_J)\right|\geq x \right] \leq  2 e^{-\frac{x^2}{2\kappa({E^i_J}+x)}}.
   \end{eqnarray*}

\subsection{Exponential bound for $\chi^2_m$} \label{ControlChi2}

We now consider $m_f$  a partition of $\{1,\dots,n\}$ such that
$\forall J \in m_f, |J|\geq \Gamma (\log(n))^2$, and we assume that
$\mathcal{M}_n$ is a set of partitions that are constructed on this
grid $m_f$. We first introduce the following set $\Omega_{m_f}$
defined by:

   \begin{eqnarray*}
     \Omf &=&\d \bigcap_{1\leq i \leq d} \Omega_{m_f}^i(\varepsilon) = \d \bigcap_{1\leq i \leq d} \bigcap_{J \in m_f} \left\{\left|\dfrac{T^i_J}{E^i_J}-1 \right|\leq \varepsilon \right\}.\\
   \end{eqnarray*}
Using the previous control, we have

\begin{eqnarray*}
     \P( \Omega_{m_f}(\varepsilon)^C)&\leq& \sum_{i=1}^d \sum_{J\in m_f}  P\left[\left|T^i_J-E^i_J\right| \geq  \varepsilon {E^i_J} \right] \leq \sum_{i=1}^d \sum_{J\in m_f}
2 e^{-\frac{\varepsilon^2{E^i_J}}{2\kappa(1+\varepsilon)}} \\&\leq&
\frac{
     C(d,\varepsilon,\Gamma,\kappa,E^{min},a)}{n^a},
   \end{eqnarray*}
with $a>2$.

The following proposition gives an exponential bound for
$\chi^2_{m}$ on the restricted event $\Omf$.
\begin{prpstn} \label{cont-chi2}
 Let $Y_1,\ldots,Y_n$ be independent random variables with distribution $s$ (from the exponential family) and verifying $(\mathbf{H})$. Let $m$ be a partition of $\mathcal{M}_n$ with $|m|$ segments and $\chi^2_{m}$ be the statistic given by (\ref{def-chi}). For any positive $x$, we have
 \begin{eqnarray*}
   \P\left[ \chi^2_m \mathbf{1}_{\Omf} \geq d\ \kappa\left( |m| + 8 (1+\varepsilon)
   \sqrt{x |{m}|} +4 (1+\varepsilon) x\right)\right] &\leq & d\ e^{-x}.
 \end{eqnarray*}
\end{prpstn}

Proof: We have
\begin{eqnarray*}
     \Eb[\chi^2(m,i)]=\sum_{J\in m}\frac{1}{{E^i_J}}\Eb(T^i_J-E^i_J)^2,
   \end{eqnarray*}
   and therefore since $ \zeta^i {E^i_J} \leq Var(T^i_J)\leq \kappa^i {E^i_J} $
   \begin{eqnarray*}
\zeta  |m| \leq \Eb[\chi^2(m,i)]\leq \kappa |m|.
   \end{eqnarray*}

We introduce the  variables $Z_J(i)$ such that
\begin{equation*}
\chi^2(m,i)=\sum_{J\in m}Z_J(i)=\sum_{J\in
m}\dfrac{(T^i_J-E^i_J)^2}{{E^i_J}},
\end{equation*}
and control their moments using
   \begin{eqnarray*}
     \E\left[ Z_J(i)^{p} \mathbf{1}_{\Omf}  \right] &\leq&\left(\frac{1}{{E^i_J}}\right)^{p} \int_{0}^{+\infty } x^{p}\ dP\left[ \left\{ (T^i_J-E^i_J)^2\geq x\right\} \cap \Omf \right] dx \\
     &\leq &\left(\frac{1}{{E^i_J}}\right)^{p} \int_{0}^{\varepsilon {E^i_J}} 2p \ x^{2p-1}P\left[|T^i_J-E^i_J|\geq x \right] dx\\
     &\leq &\left(\frac{1}{{E^i_J}}\right)^{p}  \int_{0}^{\varepsilon {E^i_J}} 4p \ x^{2p-1}e^{ -\frac{x^{2}}{2 \kappa {E^i_J} \left( 1+\varepsilon \right) }}  dx \\
     &\leq &4p\kappa^p\left( 1+\varepsilon \right) ^{p}\int_{0}^{+\infty}u^{2p-1}e^{ -\frac{u^{2}}{2}} du \\
     &\leq &4p\kappa^p\left( 1+\varepsilon \right)^{p} \int_{0}^{+\infty}\left( 2t\right) ^{p-1}e^{ -t} dt \\
     &\leq &2^{p+1}p\kappa^p\left( 1+\varepsilon\right) ^{p}p!
   \end{eqnarray*}

   We can then conclude by applying Bernstein's inequality \cite{massart2007concentration}  taking $v=2^{5}\left(\kappa( 1+\varepsilon)\right) ^{2}|m|$ and $c=4\left( \kappa(1+\varepsilon)\right)$:
   \begin{eqnarray*}
     \P\left[ \chi^2(m,i) \mathbf{1}_{\Omf} \geq \kappa\left( |m| + 8 (1+\varepsilon)
       \sqrt{x |{m}|} +4 (1+\varepsilon) x\right)\right] &\leq & e^{-x}.
   \end{eqnarray*}

While for a given $i$ the $\{Z_J(i)\}_{J\in m}$ are independent
variables, in general for a given $J$ the variables
$\{Z_J(i)\}_{1\leq i\leq d}$ are not. We conclude the proof using
lemma \ref{sum}:

\begin{lmm} \label{sum}
Let $X_1,\dots,X_n$ be real random variables and let $a_1,\dots,a_n
\in \mathbb{R}^n$ and $b_1,\dots,b_n \in [0,1]^n$  such that
$\forall 1\leq i\leq n,\ \mathbb{P}(X_i\geq a_i)\leq b_i.$ Then
$$\mathbb{P}\left(\sum_{i=1}^n X_i\geq \sum_{i=1}^n a_i\right)\leq
\sum_{i=1}^n b_i. $$
\end{lmm}

\section{Proof of our main results} \label{proof}
\subsection{Proof of proposition \ref{Risk-control}} \label{risk}

Before we focus on the main theorem, we prove our lower bound on the
risk of a model. By definition, we have

\begin{eqnarray*}
 K(\barsm,\hatsm) = \sum_{J\in m}|J|\left[A\left(\nabla A^{-1}(\bar{\T}_J)\right)-A\left(\nabla A^{-1}(\bar{\E}_J)\right)  -\bar{\E}_J . \left(\nabla A^{-1}(\bar{\T}_J)-\nabla A^{-1}(\bar{\E}_J) \right)\right].
\end{eqnarray*}

We define the compact subset $K(\varepsilon)$ of $\mathbb{R}^d$ as
the pre-image by $\nabla A$ of the domains induced by $\Omega_m$,
that is
\begin{eqnarray*}
K(\varepsilon) = \left\{\mathbf{z} \in \mathbb{R}^d \left| \nabla A(
\mathbf{z}) \in  \bigcup_{m\in\mathcal{M}}\bigcup_{J \in m}
\mathcal{B}( \bar{E}_J,\varepsilon\bar{E}_J)\right. \right\},
\end{eqnarray*}
where $\mathcal{B}(\mathbf{x},r)$ denotes the closed ball centered
in $\mathbf{x}$ with radius $r$ of $\mathbb{R}^d$. Since we consider
the union of a finite number of balls, homeomorphic properties of
$\nabla A$ ensures that $K(\varepsilon)$ is a compact set of
$\mathbb{R}^d$.\\
Since $\nabla^2 A$ is definite positive,   $A$ is
$m_{\varepsilon}$-strongly convex on the compact set
$K(\varepsilon)$, and we have \\
\begin{eqnarray*}
K(\barsm,\hatsm) \Of1 & \geq &\sum_{J\in m}|J|\dfrac{m_{\varepsilon}}{2} ||\nabla A^{-1}(\bar{\T}_J)-\nabla A^{-1}(\bar{\E}_J) ||^2 \Of1 \\
& \geq &\sum_{i=1}^d \sum_{J\in m}|J|\dfrac{m_{\varepsilon}}{2{E^i_J}} {E^i_J} \left[\nabla A^{-1}(\bar{\T}_J)^i-\nabla A^{-1}(\bar{\E}_J)^i\right]^2 \Of1, \\
\end{eqnarray*}
where $m_\varepsilon$ can be chosen as a lower bound on the smallest eigen-value of $\nabla^2 A$ on $K(\varepsilon)$. \\
Now  introducing
\begin{eqnarray} \label{def-V2}
V^2_m =\sum_{i=1}^d V^2(m,i) = \sum_{i=1}^d\sum_{J\in
m}{E^i_J}\left[\nabla A^{-1}(\bar{\E}_J)^i-\nabla
A^{-1}(\bar{\Tb}_J)^i\right]^2,
\end{eqnarray}
we obtain
\begin{eqnarray} \label{K2-V2}
K(\barsm,\hatsm) \Of1 \geq \dfrac{m_{\varepsilon}}{2E^{max}} V^2_m
\Of1.
\end{eqnarray}

{Now if we consider $M_{\varepsilon}$  an upper bound of the
eigen-values of $\nabla^2 A$ on $K(\varepsilon)$, then on $\Omf$ we have}

$$||\nabla A^{-1}(\bar{\T}_J)-\nabla A^{-1}(\bar{\E}_J) ||^2 \geq
\frac{1}{M_{\varepsilon}} ||\bar{\T}_J-\bar{\E}_J ||^2,$$

and therefore
\begin{eqnarray*}
K(\barsm,\hatsm) \Of1 & \geq &\sum_{J\in m}|J|\dfrac{m_{\varepsilon}}{2M_{\varepsilon}} ||\bar{\T}_J-\bar{\E}_J ||^2 \Of1 \\
& \geq &\dfrac{E^{min} m_{\varepsilon}}{2M_{\varepsilon}}\chi^2_{m} \Of1.
\end{eqnarray*}


Combining the previous results leads to
\begin{eqnarray*}
\frac{1}{2}E^{min}\dfrac{ m_{\varepsilon}}{M_\varepsilon}\
\mathbb{E}[\chi^2_m]-\mathbb{E}[\chi^2_m\mathbf{1}_{\Omf^c}]\leq\mathbb{E}[K(\barsm,\hatsm)\Of1],
\end{eqnarray*}

and using $\mathbb{E}[\chi^2_m]\geq d\zeta |m|$ and Cauchy-Schwarz
inequality on $\mathbb{E}[\chi^2_m\mathbf{1}_{\Omf^c}]$,

\begin{eqnarray*}
\frac{1}{2}\xi E^{min}\dfrac{ m_{\varepsilon}}{M_\varepsilon}\ d
|m|-\dfrac{C(a,\Gamma,E^{min},E^{max},\varepsilon,\kappa,d)}{n^{a/2-1}}\leq\mathbb{E}[K(\barsm,\hatsm)\Of1].
\end{eqnarray*}
Finally, introducing $C(\varepsilon)=\frac{1}{2}\xi E^{min}\dfrac{ m_{\varepsilon}}{M_\varepsilon}$, and since
$\mathbb{E}\left[K(\barsm,\hatsm)\mathbf{1}_{\Omf^C} \right] \geq
0,$ we have
\begin{eqnarray*}
K(s,\barsm) + C(\varepsilon)\ d
|m|-\dfrac{C(a,\Gamma,E^{min},E^{max},\varepsilon,\kappa,d)}{n^{a/2-1}}\leq\mathbb{E}[K(s,\hatsm)].
\end{eqnarray*}

\subsection{Proof of theorem \ref{theo}}

To prove theorem \ref{theo}, we will control each of the three terms
in equation (\ref{decomposition}) individually. To this effect, we
introduce a generic partition $m'$ of $\mathcal{M}_n$.

\subsubsection{Control of term   $\barga(\hatsmp)-\barga(\barsmp) $}

Let
\begin{eqnarray*}
 \OO1 = \bigcap_{m' \in \M_n} \left\{ \chi^2_{m'}  \mathbf{1}_{\Omf} \right. \leq d \kappa  \left[ |m'| \right. &+& 8(1+\varepsilon) \sqrt{(L_{m'} |m'|+\xi)|m'|}  \\
   &+ &  \left. \left. 4(1+\varepsilon)(L_{m'}|m'|+\xi)\right]\right\}.
\end{eqnarray*}
 Then we can
control   $\barga(\hatsmp)-\barga(\barsmp) $ with the following
proposition :
\begin{prpstn}  \label{t3}
     \begin{eqnarray*}
       &&\! \left(\bar{\gamma}(\hat{s}_{m'})-\bar{\gamma}(\bar{s}_{m'})\right) \mathbf{1}_{\Omf\cap\OO1}  \leq \frac{1}{1+\varepsilon}K(\bar{s}_{m'},\hat{s}_{m'}) + \\
       &&+ \dfrac{dE^{max} \kappa(1+\varepsilon)}{2m_{\varepsilon}} \left[|m'|  +8(1+\varepsilon) \sqrt{(L_{m'} |m'|+\xi)|m'|} +  4(1+\varepsilon)(L_{m'}|m'|+\xi)  \right]
     \end{eqnarray*}
\end{prpstn}
Proof: We recall that
\begin{eqnarray*}
\barga(\barsmp)&=&\sum_{J\in m'}|J|\left[-\nabla A^{-1}(\bar{\E}_J).(\bar{\Tb}_J-\bar{\E}_J) \right] \text{and}\\
\barga(\hatsmp)&=&\sum_{J\in m'}|J|\left[-\nabla
A^{-1}(\bar{\Tb}_J).(\bar{\Tb}_J-\bar{\E}_J) \right].
\end{eqnarray*}
We therefore have
\begin{eqnarray*}
\barga(\barsmp)-\barga(\hatsmp)&=&\sum_{J\in m'}|J|\left[(\bar{\E}_J-\bar{\Tb}_J).\left(\nabla A^{-1}(\barT_J)-\nabla A^{-1}(\barE_J)\right) \right]\\
&=&\sum_{i=1}^d\sum_{J\in
m'}|J|\left[(\bar{E^i}_J-\bar{T^i}_J)\left(\nabla
A^{-1}(\barT_J)^i-\nabla A^{-1}(\barE_J)^i\right)\right].
\end{eqnarray*}

Using Cauchy-Schwarz inequality,
\begin{eqnarray*}
\barga(\barsmp)-\barga(\hatsmp)&=&\sum_{i=1}^d\sum_{J\in m'}\left[\dfrac{({E^i}_J-{T^i}_J)}{\sqrt{{E^i_J}}}\sqrt{{E^i_J}}\left(\nabla A^{-1}(\barT_J)^i-\nabla A^{-1}(\barE_J)^i\right)\right]\\
&\leq & \left[\sqrt{\sum_{i=1}^d\sum_{J\in m'}\dfrac{({T^i}_J-{E^i}_J)^2}{{E^i_J}}} \sqrt{\sum_{i=1}^d\sum_{J\in m'}{E^i_J}\left[\nabla A^{-1}(\barT_J)^i-\nabla A^{-1}(\barE_J)^i \right]^2} \ \right] \\
&\leq &  \sqrt{\chi^2_{m'}} \sqrt{V^2_{m'}},
\end{eqnarray*}
where $\chi^2_m$ has been defined in equation (\ref{def-chi}) and
$V^2_m$ has been defined in equation (\ref{def-V2}). Then from
equation (\ref{K2-V2}) and $2ab\leq xa^2+x^{-1}b^2$, we get on
$\Omf$

\begin{eqnarray*}
\barga(\barsmp)-\barga(\hatsmp)&\leq&
\dfrac{1+\varepsilon}{2m_{\varepsilon}}E^{max}\chi^2_{m'} +
\dfrac{1}{1+\varepsilon}K(\bar{s}_{m'},\hatsmp).
\end{eqnarray*}
Finally, using proposition \ref{cont-chi2},
     \begin{eqnarray*}
       \left(\bar{\gamma}(\hat{s}_{m'})-\bar{\gamma}(\bar{s}_{m'})\right) \mathbf{1}_{\Omf\cap\OO1} & \leq &  \dfrac{d \kappa(1+\varepsilon)}{2m_{\varepsilon}}E^{max} \left[|m'| +8(1+\varepsilon) \sqrt{(L_{m'} |m'|+\xi)|m'|} \right.\\
       & & + \left. 4(1+\varepsilon)(L_{m'}|m'|+\xi)  \right] +
       \frac{1}{1+\varepsilon}K(\bar{s}_{m'},\hat{s}_{m'}).
     \end{eqnarray*}

\subsubsection{Control of term   $\barga(\barsm)-\barga(s) $}
The expectation of the second term can be bounded using the
following proposition:
\begin{prpstn}  \label{t2}
\begin{eqnarray*}
       \left|\mathbb{E}\left[(\barga(\barsm)-\barga(s))\mathbf{1}_{\Omf}  \right]\right| &\leq&
       \dfrac{C(d,\varepsilon,\Gamma,\kappa,\nu,E^{min},E^{\max})}{n^{(a-1)/2}},
     \end{eqnarray*}
where $\nu$ is the greatest length of the compact $\Nu$.
\end{prpstn}
Proof: We recall that
\begin{eqnarray*}
\barga(\barsm)-\barga(s) = \sum_{J\in m}\sum_{t\in J}
(\thetabf_t-\nabla A^{-1}(\barE_J)). (\Tb_t-\E_t).
\end{eqnarray*}
Then
\begin{eqnarray*}
\left|\Eb\left[(\barga(\barsm)-\barga(s))\mathbf{1}_{\Omf}\right]\right| &\leq &\Eb\left[\left|\barga(\barsm)-\barga(s) \right|\mathbf{1}_{\Omf}^C\right] \\
&\leq& \Eb\left[ \left|\sum_{J\in m}\sum_{t\in J} \sum_{i=1}^d (\theta^i_t-\nabla A^{-1}(\barE_J)^i) (T^i_t-E^i_t)\right|\mathbf{1}_{{\Omf}^C}\right]\\
&\leq& \nu \sum_{i=1}^d \left[\Eb\left(\sum_{t=1}^n(T^i_t-E^i_t) \right)^2 \right]^{1/2}\mathbb{P}\left(\Omega_{m_f}^i (\varepsilon)^C\right)^{1/2}\\
&\leq& \nu \sum_{i=1}^d \left[\sum_{t=1}^n \kappa {E^i_t} \right]^{1/2}\mathbb{P}\left(\Omega_{m_f}^i (\varepsilon)^C\right)^{1/2}\\
&\leq&  d\nu\dfrac{\sqrt{2 \kappa
E^{\max}C(\varepsilon,\Gamma,\kappa,E^{min})}}{n^{(a-1)/2}}.
\end{eqnarray*}

\subsubsection{Control of term   $\barga(\hatsmp)-\barga(\barsmp) $}
Finally, we control the last term using proposition \ref{s-u} for
which a proof can be found in \cite{CleynenLebarbier2014}.

\begin{prpstn} \label{s-u}
       \begin{eqnarray*}
       \P\left[\bar{\gamma}(s)-\bar{\gamma}(u) \geq  K(s,u) - 2 h^2(s,u)+2dx \right] &\leq & e^{-dx} \leq d e^{-x}.
     \end{eqnarray*}
\end{prpstn}

Applying it to $u=\bar{s}_{m'}$ yields:
   \begin{eqnarray*} 
     \P\left[\bar{\gamma}(s)-\bar{\gamma}(\bar{s}_{m'})\geq K(s,\bar{s}_{m'})-2 h^2(s,\bar{s}_{m'})+2dx \right] &\leq &  d e^{-x}.
   \end{eqnarray*}

   We then define $$\mO2=  \bigcap_{m' \in \M_n}\left\{ \bar{\gamma}(s)-\bar{\gamma}(\bar{s}_{m'}) \leq   K(s,\bar{s}_{m'}) - 2  h^2(s,\bar{s}_{m'})+2d(L_{m'}|m'|+\xi) \right\}.$$

\subsubsection{Proof of the theorem}\label{prooftheo}

We can now combine the previous propositions in order to prove our
main theorem. To this effect, we introduce the following event:
$\Omega(\varepsilon,\xi) = \Omf\cap\OO1 \cap \mO2$. On this event,
we have, with $R= \bar{\gamma}(\bar{s}_m)-\bar{\gamma}(s)$:
   \begin{eqnarray*}
     (\barga(\barsm) -\barga(\hatsmh))\O & = & (\barga(s)-\barga(\barsmh))\O +(\barga(\barsm) -\barga(s))\O +(\barga(\barsmh)-\barga(\hatsmh))\O\\
&\leq & \left[ K(s,\barsmh)-2 h^2(s,\barsmh)\right]\O + R\O + \frac{1 }{1+\varepsilon}K(\barsmh,\hatsmh)\O \\
     &  & +  \dfrac{dE^{max} \kappa(1+\varepsilon)}{2m_{\varepsilon}}\left[|\mh| +8(1+\varepsilon) \sqrt{(L_{\mh} |\mh|+\xi)|\mh|}+ 4(1+\varepsilon)(L_{\mh}|\mh|+\xi) \right] \\
     &&+ 2d L_{\mh} |\mh| +2d\xi.
   \end{eqnarray*}

We then apply the  proof of  \cite{CleynenLebarbier2014} with
$C(\varepsilon)=
\dfrac{E^{max}\kappa(1+\varepsilon)}{2m_{\varepsilon}}$, which
yields:
   \begin{eqnarray*}
      \dfrac{\varepsilon}{1+\varepsilon} h^2(s,\hatsmh)\O &\leq&  K(s,\barsm)\O + R\O -pen(\hat{m})+pen(m)\\
     & & + d |\mh| C_2(\varepsilon) \left(1 + 4\sqrt{L_{\mh}} \right)^2 + 2\xi d\left[1
     +(1+\varepsilon)C(\varepsilon)\left(\dfrac{8}{\varepsilon}+2\right)\right],
   \end{eqnarray*}
where $C_2(\varepsilon) =
\dfrac{E^{max}\kappa}{2m_{\varepsilon}}(1+\varepsilon)^3$. Then
since by assumption, $pen(\mh)\geq \beta_{\varepsilon} d|\mh|
\left(1 + 4\sqrt{L_{\mh}} \right)^2$, choosing  $\beta_{\varepsilon}
= C_2(\varepsilon)$ yields:
   \begin{eqnarray*}
      h^2(s,\hat{s}_{\hat{m}})\O &\leq& C_{\beta_{\varepsilon}} \left[ K(s,\barsm)\O + R\O + pen(m)\right]+d \xi C(\beta_{\varepsilon}),
   \end{eqnarray*}
with $C_{\beta_{\varepsilon}} =\dfrac{1+C_2^{-1}(\beta_\varepsilon)}{C_2^{-1}(\beta_\varepsilon)}$. \\

From propositions \ref{cont-chi2} and \ref{s-u} comes both
$\P\left(\OO1^C \right) \leq d\sum_{m' \in \M_n} e^{-(L_{m'}|m'|
+\xi)}$  and $\P\left(\mO2^C \right) \leq d\sum_{m' \in \M_n}
e^{-(L_{m'}|m'| + \xi)}$, so that using hypothesis (\ref{weights}),
   \begin{eqnarray*}
     \P\left(\OO1^C \cup \mO2^C \right) &\leq& 2 d\sum_{m' \in \M_n} e^{-(L_{m'}|m'| + \xi)} \leq 2d\Sigma e^{-\xi},
   \end{eqnarray*}
   and thus $ \P\left(\OO1 \cap \mO2 \right) \geq 1- 2d\Sigma e^{-d\xi}$. Integrating over $\xi$ and using proposition \ref{t2} leads to

\begin{eqnarray*}
     \E\left[ h^2(s,\hatsmh)\mathbf{1}_{\Omf} \right] &\leq & C_{\beta_{\varepsilon}}  \left[ K(s,\barsm)+ \dfrac{C(d,\nu,\beta_{\varepsilon},\Gamma,\kappa,E^{min},E^{\max},a)}{n^{(a-1)/2}}+pen(m)\right] + 2d\Sigma C(\beta_{\varepsilon}).
   \end{eqnarray*}

   And since  $ \E\left[ h^2(s,\hatsmh)\mathbf{1}_{\Omf^C} \right] \leq \dfrac{C(d,\beta_{\varepsilon},\Gamma,\kappa,E^{min},a)}{n^{a-1}}$, we have

   \begin{eqnarray*}
     \E\left[ h^2(s,\hatsmh) \right] &\leq & C_{\beta_{\varepsilon}} \left[ K(s,\barsm) +pen(m)\right] + C'(d,\nu,\beta_{\varepsilon},\Gamma,\kappa,E^{min},E^{\max},\Sigma).
   \end{eqnarray*}

   Finally, by minimizing over $m\in\mathcal{M}_n$, we get
       \begin{eqnarray*}
     \E\left[ h^2(s,\hatsmh) \right] &\leq & C_{\beta_{\varepsilon}} \inf_{m\in\mathcal{M}_n} \left\{ K(s,\barsm) +pen(m)\right\} +  C'(d,\nu,\beta_{\varepsilon},\Gamma,\kappa,E^{min},E^{\max},\Sigma).
   \end{eqnarray*}

\section{Characteristics of classic laws}\label{classic}

In this section we provide some explicit values for the constants used in the main theorem for some usual distributions from the exponential family. While we do not detail all computations, in Table \ref{usual} we summarize the parameters of the distribution, the sufficient statistics associated with the natural parameters, the expression of the log-partition function and possible choices for the values of $\zeta$ and $\kappa$. \\

Here we detail the computations for the exponential distribution. Some details for other distributions are provided in Appendix \ref{app}. \\
The exponential distribution can be written in its natural form as $$s(t)=\exp\left(-\lambda_t Y_t +\log \lambda_t \right) $$ from which we obtain, dropping the index $t$
\begin{eqnarray*}
  \theta &=& -\lambda \\
  A(\theta)&=&-\log(-\theta)\\
  T(Y)&=&Y.
\end{eqnarray*}

\begin{itemize}
\item Variance / Expectation relationship \\
$A'(\theta)=-\frac{1}{\theta} $ and $A"(\theta)=\frac{1}{\theta^2}$  so that $$-\frac{1}{\theta^{min}} E \leq Var(T)\leq -\frac{1}{\theta^{max}} E .$$

\item Cumulants \\
$A(\theta+z)-A(\theta) =-\log \left(1+\frac{z}{\theta} \right)$ is analytic for $-|\theta|<z<|\theta|$, \textit{i.e.} $R=|\theta|$. We then obtain
$$A(\theta+z)-A(\theta) = \sum_{k\geq 1}(-1)^k\frac{z^k}{k\theta^k} = \sum_{k\geq 1}\left(-\frac{1}{\theta}\right)^k(k-1)!\frac{z^k}{k!} $$ so that for $k\geq 1, c_k=\left(-\frac{1}{\theta}\right)^k(k-1)!.$
\item Bounds \\
We have, for $k\geq 3$;
$$\left|\frac{c_k}{c_2^{k/2}} \right|=(k-1)! \leq \frac{k}{2} (k-1)! \leq \frac{1}{2} k!(\alpha)^{k-2} $$
with $\alpha=1$. Then $\tau=\max\left\{ \alpha \sqrt{Var(T)},\frac{1}{R} \right\} = \max \left\{ \frac{1}{|\theta|}, \frac{1}{|\theta|} \right\} = -\frac{1}{\theta^{max}}$ and finally, $\zeta=-\frac{1}{\theta^{min}} $ and $\kappa=-\frac{1}{\theta^{max}}$.
\end{itemize}

\begin{table}[h!]
\begin{center}
{\footnotesize
\begin{tabular}{M{1.6cm}|M{1.6cm}M{1.3cm}M{2.8cm}M{2.8cm}M{2.8cm}N}
Distribution &  Natural Parameters  &  Sufficient Statistics& Log-Partition & $\zeta$ & $\kappa$ &\\
\hline
Poisson & $\log \lambda$  & $x$ & $e^\theta$ & $1$ & $1$ &\\[20pt]
\hline
Exponential  & $-\lambda$  & $x$ & $-\log(-\theta)$ & $-\dfrac{1}{\theta^{min}}$ & $-\dfrac{1}{\theta^{max}}$ &\\[20pt]
\hline
Gaussian \small{($\mu>0$)}  & $\dfrac{\mu}{\sigma^2},-\dfrac{1}{2\sigma^2}$  & $(x,x^2)$ &$-\dfrac{{\theta_1}^2}{4\theta_2}-\dfrac{1}{2}\log(-2\theta_2)$ &  $\min\left\{\frac{1}{\theta_1^{max}},-\frac{1}{\theta_2^{min}}\right\}$ & $\max\left\{\frac{1}{\theta_1^{min}},-\frac{2}{\theta_2^{max}}\right\}$ &\\[20pt]
\hline
Pareto \small{($\alpha<\frac{-1}{log x_m}$)} & $-\alpha$ & $\log x$ & $\theta \log x_m - \log(-\theta)$ &$\dfrac{1}{\theta^{max}}\dfrac{1}{{\theta^{max}}\log x_m-1}$ & $\frac{\log x_m}{\theta^{min}\log x_m-1}$& \\[20pt]
\hline
Gamma \hspace{1cm} \small{($\alpha$ fixed)} & $-\beta $ & $x$ & -$\alpha\log(-\theta)$ & $-\dfrac{1}{\theta^{min}} $& $-\dfrac{1}{\theta^{max}} $ & \\[20pt]
\hline
Weibull \hspace{1cm} \small{($k$ fixed)} & $-\lambda^{-k}$  & $x^k$ & $-\log(-\theta)$ & $-\dfrac{1}{\theta^{min}}$ & $-\dfrac{1}{\theta^{max}}$ &\\[20pt]
\hline
Laplace \hspace{1cm} \small{($\mu$ fixed)} & $-1/b$  & $|x-\mu|$ & $-\log(-\theta)$ & $-\dfrac{1}{\theta^{min}}$ & $-\dfrac{1}{\theta^{max}}$ &\\[20pt]
\hline
Binomial \hspace{1cm} \small{($n$ fixed)} & $\log{(\frac{p}{1-p})}$  & $x$ & $n\log(1+e^\theta)$ & $\dfrac{1}{1+e^{\theta_{max}}}$ & $2$ &\\[20pt]
\hline
Negative Binomial \small{($\phi$ fixed)} & $\log{(1-p)}$  & $x$ & $-\phi\log(1-e^\theta)$ & $1$ & $e^{-\theta_{min}}$ &\\[20pt]
\hline
\end{tabular}
}
\end{center}
\caption{List of most usual univariate laws from the exponential
family.}\label{usual}
\end{table}

\section{Particular case of categorical variables}\label{categorical}
The goal of this section is double. We first wish to illustrate how
our general approach can be used when working with a specific
distribution, in this  example the categorical distribution. Indeed,
while we have shown that the penalty shape is the same for all
distribution, the constants can be defined properly with quantities
depending on the distribution of interest. We then show how the main
results could be derived if working directly with the
characteristics of the categorical distribution instead of dealing
with the log-partition function from the exponential family. We
conclude this section by comparing the constants obtained in each
case.

\subsection{Application of the general approach} \label{cat:general}

Here we suppose that $Y$ can take values between $1$ and $d+1$ and
we denote $\{p^i;1\leq i\leq d\}$ the probability that $Y$ belongs
to  categories $1$ through $d$ (so that $p_{d+1}=1-\sum_{i=1}^d
p^i$).\newline

In the canonical form, the parameters $\boldsymbol{\theta}$ are
given by $\theta^i=\log\frac{p^i}{1-\sum_{i=1}^d p^i}$ (for $1\leq
i\leq d$) and we have
\begin{itemize}
\item $T^i=\mathbf{1}_{\{Y=i\}}$
\item $A(\thetabf)=\log(1+\sum_{i=1}^d e^{\theta^i}) $
\item $E^i=\dfrac{e^{\theta^i}}{1+\sum_{j=1}^d e^{\theta^j}}$ and
\begin{eqnarray*}
 \nabla^2 A(\thetabf)=\left(
   \begin{array}{ccccc}
     E^1(1 -E^1) & -E^1E^2 & &\dots&-E^1E^d\\
     -E^1E^2 &E^2(1 -E^2)  & &\dots &-E^2E^d\\
     \vdots&&&&\\ &&&&\\
     -E^1E^d &\dots  & & &E^d(1-E^d)\\
   \end{array} \right)
\end{eqnarray*}
\end{itemize}
The inverse mapping of the gradient of the log-partition function is
given by $[(\nabla A)^{-1}(\boldsymbol{\beta})]^i =\log
\dfrac{\beta^i}{1-\sum_{j=1}^d \beta^j}.$

\subsubsection{Verification of hypothesis (H)}
Let us assume that the probabilities of categories  $1$ to $d$ are
bounded away from $0$ and $1$. We will want to assume the same
property on the \textit{normalizing}  category $d+1$, so that there
exists $a>1$ such that $-\log a < \theta_i < \log a$ for all $1\leq
i \leq d$. Our compact set $\boldsymbol{\nu}$ is therefore defined
by $[-\log a;\log a]^d$, and $\nu=2\log a$. This leads to:
\begin{eqnarray*}
E^{min}=\dfrac{1}{a+a^2d}>0 \qquad \text{and} \qquad
E^{max}=\dfrac{a^2}{d+a}<1.
\end{eqnarray*}

Now since for all $i$ we have $Var(T^i)=(1-E^i)E^i$, we get
$\dfrac{d+a(1-a)}{d+a}E^i \leq VarT^i \leq E^i $ and therefore
$$\zeta= \dfrac{d+a(1-a)}{d+a}.$$

The Laplace transform of sufficient statistic $T^i$ is
\begin{eqnarray*}
 \log \mathbb{E}\left[e^{z(T^i-E^i)} \right]&=&\log \left [\dfrac{1+\sum_j e^{\theta^j}+e^{\theta^i}\left(e^z-1 \right)}{1+\sum_j e^{\theta^j}}\right ]
 = \log \left[1 + E^i\left(e^z-1\right)\right],
\end{eqnarray*}
which is analytic in $z$ provided $z\leq \log\left(1+\frac{1}{E^i}
\right).$

Let $R=\log 2$ and let $z\in[-R;R]$. Then the  cumulants associated
with $T^i$ can be obtained through the recurrence property
$c_1^i=p^i$ and for $k\geq 2$,
\begin{eqnarray*}
c_{k+1}^i=p^i(1-p^i)\dfrac{\partial c_k^i}{\partial p^i},
\end{eqnarray*}

where we have switched back to usual proportion notations for sake of readability.\\
Let {$P_{k+1}(p)=\frac{\partial c_k^i}{\partial p}$} with $p=p^i$.
We have $P_2(p)=1$, polynomial in $p$ of degree $0=k-2$ and with sum
of absolute coefficients $=1=2^{2-3}*2!$. Let us assume that for a
given $k>1$, $P_k(p)$ is a polynomial in $p$ of degree $k-2$ and
with sum of absolute coefficients less than $2^{k-3}k!$. Then
\begin{eqnarray*}
P_{k+1}(p)=\frac{\partial c_{k}^i}{\partial p}=\frac{\partial
\left[p(1-p)P_k(p) \right]}{\partial
p}=(1-2p)P_k(p)+(p-p^2)\frac{\partial P_k(p)}{\partial p}
\end{eqnarray*}

Thus denoting $P_k(p)=\sum_{\ell=0}^{k-2}a_\ell p^\ell$ we get
$P_{k+1}(p)=a_0+\sum_{\ell=1}^{k-2}(\ell+1)(a_\ell-a_{\ell-1})p^\ell-ka_{k-2}p^{k-1}$
which is of degree $k-1$ and with sum of absolute coefficient less
than

\begin{eqnarray*}
|a_0|+\sum_{\ell=1}^{k-2}(\ell+1)|a_\ell-a_{\ell-1}|+k|a_{k-2}| &\leq & 3|a_0|+ \sum_{\ell=1}^{k-3}(2*\ell+3)|a_\ell|+(2*k-1)|a_{k-2}| \\
&\leq &\sum_{\ell=0}^{k-2}(2*\ell+3)|a_\ell| \leq  2*(k+1)
\sum_{\ell=0}^{k-2}|a_\ell| \leq 2^{k-2}(k+1)!
\end{eqnarray*}

We therefore obtain by recurrence that $P_k$ (for $k\geq 2$) is a
polynomial in $p$ of degree $k-2$ and with sum of absolute
coefficients less than $2^{k-3}k!$. From this we get:
\begin{eqnarray*}
 \left|\dfrac{c_k(T^i)}{c_2(T^i)^{k/2}}\right| \leq \frac{1}{2}k!\left(\frac{2}{\sqrt{p^i(1-p^i)}}\right)^{k-2}
\end{eqnarray*}
and
$\tau_i=\max\left\{\frac{2}{\sqrt{p^i(1-p^i)}}*\sqrt{p^i(1-p^i)},\frac{1}{R^i}\right\}=\max\{2,1/R\}.$
Since $\kappa_i=1$ and $R=\log 2$, we can take $\kappa = 2$.

Note that considering the sufficient statistic $T^i$ independently
of other categories can be interpreted as considering a surrogate
random variable $X$ with success $X=1$ if $Y=1$ and $X=0$ otherwise.
The distribution of $X$ is then Bernoulli with parameter $p=p^i$ and
the results claimed in the previous Section can be deducted from
what precedes.

\subsubsection{Computation of $m_{\varepsilon}$ and  $M_{\varepsilon}$}

For a given vector of proportions $(\pi_1,\dots,\pi_d)$, the
smallest eigen value of $\nabla^2 A$ is greater than $\Pi
=(1-\sum_{i=1}^d \pi_i) \prod_{i=1}^d \pi_i$. This can be seen by
applying Sylvester's criterion to the symmetric matrix $\nabla^2 A-
\Pi Id$. Considering the mapping between natural parameters and
proportions, we obtain that $\Pi$ can be defined as $\Pi(\eta) =
\dfrac{\exp(\sum \eta^i)}{\left[1+\sum e^{\eta^i}
\right]^{d+1}}.$\newline

Therefore let $z \in K(\varepsilon)$. There exists $E$ such that
$\nabla A(z) \in \mathcal{B}(E,\varepsilon)$, that is there exists
$E$ such that $\forall i, \frac{e^{z^i}}{1+\sum e^{j^j}} \in
[E^i-\varepsilon; E^i+\varepsilon]$. We therefore have:
\begin{eqnarray*}
\Pi(z)=\dfrac{\prod e^{z^i}}{\left[1+\sum e^{z^i} \right]^{d+1}}
\geq (E^{min}-\varepsilon)^{d+1}
\end{eqnarray*}
and can conclude by taking
$m_{\varepsilon}=(E^{min}-\varepsilon)^{d+1}.$ \newline

The greatest eigen-value of $\nabla^2 A$ can be bounded by the trace of the gradient matrix, in this case the sum of the variance of each sufficient statistic. Therefore, for a distribution with proportions $\{q_i\}$, one would get
\begin{eqnarray*}
M\leq \sum_{i=1}^d q_i(1-q_i)\leq \sum_{i=1}^d q_i \leq 1.
\end{eqnarray*}
Therefore, on any set we have $M\leq 1$ and in particular, on $K(\varepsilon)$, we can take $M_\varepsilon=1$.

\subsection{Direct results for the categorical variables}\label{categoricalDirect}

This section is dedicated to the derivation of direct controls when
studying the categorical distribution. As before, $Y$ will take
values in $1,\dots, d+1$ and $p^i$ will denote the probability of
category $i$. Once again, it is possible to  reduce the number of
parameters to $d$, however  keeping all $d+1$ parameters leads to
more tractable quantities to control and  to smaller resulting
constants. For sake of readability and to avoid confusions, we will
denote $r=d+1$.

\subsubsection{Notations and main result}
The density $s_t$ of $Y_t$ can be decomposed in $r$ terms by
denoting $s(t,i)=P\left(
Y_{t}=i\right)=p^i_t=E[\mathbf{1}_{\{Y=i\}}]$. In this specific
case,  quantities such as the contrast or the Kullback-Leibler risk
can be identified :
\begin{itemize}
\item The model $S_m$ is:
\begin{equation*}
 \mathcal{S}_{m}=\left\{\begin{array}{c}
 u:\left\{ 1,...,n\right\} \times \left\{ 1,...,r\right\} \rightarrow \left[0,1\right] \text{ such that} \\
 \forall J\in m\text{, }\forall i\in \left\{ 1,2,...,r\right\},u\left( t,i\right) =u\left( t^{\prime },i\right) \text{ }=u\left(
 J,i\right) \text{\ \ }\forall \text{ }t,t^{\prime }\in J
 \end{array} \right\} . \nonumber
\end{equation*}
\item the contrast is the $\log$-likelihood defined by:
\begin{eqnarray*}
 \gamma\left( s\right) &=&-\log \left[ \prod_{t=1}^{n}\prod_{i=1}^{r} s\left( t,i\right) ^{\mathbf{1}_{\{Y=i\}} }\right]=-\sum_{t=1}^{n}\sum_{i=1}^{r} \mathbf{1}_{\{Y=i\}} \log \left[ s\left( t,i\right) \right] \text{.}
\end{eqnarray*}
\item Associated to this contrast, the Kullback-Leibler information between $s$ and $u$ is:
\begin{equation*}
K\left( s,u\right) =\sum_{t=1}^{n}\sum_{i=1}^{r}s\left( t,i\right)
\log \left[ \frac{s\left( t,i\right) }{u\left( t,i\right) }\right]
\end{equation*}
\item For a given partition $m$, we obtain the minimum contrast estimator of $s$
\begin{equation*}
\hat{s}_{m}\left( t,i\right) =\frac{N_{J}\left( i\right) }{|\text{
}J\text{ }|},\text{ for }i\in \left\{ 1,...,r\right\} ,t\in J\text{
and }J\in m\text{,} \nonumber
\end{equation*}%
with $N_{J}\left( i\right) =\sum_{t\in J} \mathbf{1}_{\{Y=i\}} $,
and the projection
\begin{eqnarray*}
\bar{s}_{m}\left( t,i\right) &=& \arg\min_{u \in \mathcal{S}_{m}}
K(s,u) =\frac{\sum_{t\in J}s\left( t,i\right) }{\left\vert \text{
}J\text{ }\right\vert }  \nonumber
\end{eqnarray*}
\end{itemize}

The following theorem gives the direct version of theorem \ref{theo}.
\begin{thrm} \label{theo_categorial_direct}
Suppose that one observes independent variables $Y_{1}$, ..., $Y_{n}$ taking their values in $\left\{1,2,...,r\right\} $ with $r\in \mathbb{N}$ and $r\geq 2$. We define for $t\in \left\{ 1,...n\right\} $ and $i\in \left\{ 1,2,...,r\right\} $
\begin{equation*}
P\left( Y_{t}=i\right) =s\left( t,i\right)
\end{equation*}
and consider a collection $\mathcal{M}_{n}\mathcal{\ }$ of partitions constructed on the grid $\left\{ 1,...,n\right\} $. Let
$\left( L_{m}\right)_{m\in \mathcal{M}_{n}}$ be some family of positive weights and define $\Sigma $ as
\begin{equation*}
\Sigma =\sum_{m\in \mathcal{M}_{n}}\exp \left( -L_{m}|m|\right)
<+\infty \text{.}  \nonumber
\end{equation*}
Assume that
\begin{itemize}
 \item there exists some positive absolute constant $\rho $ such that $s\geq \rho $,
 \item $\mathcal{M}_{n}$ is a collection of partitions constructed on a partition $m_{f}$ such that $|J|\geq \Gamma \left[ \log \left( n\right) \right] ^{2}\forall J\in m_{f}$ where $\Gamma $ is a positive absolute constant.
\end{itemize}
Let $\lambda_{\varepsilon} >1/2$. If for every $m\in
\mathcal{M}_{n}$
\begin{equation} \label{pen_generale_likelihood}
 {pen}\left( m\right) \geq \lambda_{\varepsilon}  \text{ }r|m|\text{ }\left( 1+4\text{ }\sqrt{L_{m}}\right)^2,
\end{equation}
then
\begin{equation} \label{risk-direct}
 \Eb\left[ h^{2}\left( s,\hatsmh\right) \right] \leq C_{\lambda_{\varepsilon} }\inf_{m\in \mathcal{M}_{n}}\left\{ K\left( s,\overline{s}_{m}\right) +{pen}\left(m\right) \right\} +C\left( \Sigma ,r,\lambda_{\varepsilon}  ,\rho ,\Gamma \right) ,
\end{equation}%
with $C_{\lambda_{\varepsilon}
}=\frac{2\text{{}}(2\lambda_{\varepsilon} )
^{1/3}}{(2\lambda_{\varepsilon} ) ^{1/3}\text{{}}-1}$.
\end{thrm}

We obtain this result by following the same lines of the proof of
theorem \ref{theo}. The main constant of interest is
$C_{\lambda_{\varepsilon} }$ which comes from the control of the
term $\bar{\gamma}\left( \bar{s}_{{m}^{\prime
}}\right)-\bar{\gamma}\left( \hat{s}_{m^{\prime }}\right)$ on a
particular set $\Omf$. Here again this control is obtained using two
results: the control of a chi-square statistic $\chi _{m}^{2}$ and
the control $K\left(\bar{s}_{m},\hat{s}_{m}\right)$ through a
quantity $V_{m}^{2}$ on  $\Omf$.

\subsubsection{Control of $\bar{\gamma}\left( \bar{s}_{{m}^{\prime }}\right)-\bar{\gamma}\left( \hat{s}_{m^{\prime }}\right)$}

We write for
$m^{\prime}\in \mathcal{M}_{n}$%
\begin{equation*}
 \bar{\gamma}\left( \bar{s}_{m^{\prime }}\right) -\bar{\gamma}\left( \hat{s}_{m^{\prime }}\right)
 =\sum_{t=1}^{n}\sum_{i=1}^{r}\overline{\mathbf{1}_{\left\{ Y_{t}=i\right\}}}\log \left[ \frac{\hat{s}_{m^{\prime }}\left( t,i\right)}{\bar{s}_{m^{\prime }}\left( t,i\right) }\right] \text{.}
\end{equation*}%
By Cauchy-Schwarz inequality,
\begin{eqnarray*}
 \bar{\gamma}\left( \bar{s}_{m^{\prime }}\right) -\bar{\gamma}\left(\hatsmp\right) &\leq &\sqrt{\sum_{i=1}^{r}\chi _{m^{\prime }}^{2}\left( i\right)}\times V_{m^{\prime }}^2\text{,} \notag
\end{eqnarray*}
where
\begin{equation*}
 \chi _{m}^{2}=\sum_{i=1}^r \chi _{(m,i)}^{2} =\sum_{i=1}^r  \sum_{J\in m} Z_J (i),
\end{equation*}
with $Z_J (i)=\frac{\overline{N_{J}\left( i\right)}^{2}}{\sum_{t\in
J}s\left( t,i\right) }$ and $\overline{N_{J}\left( i\right)
}=\sum_{t\in J}\ \left[\mathbf{1}_{\{Y=i\}}-s(t,i)\right]$,  and
\begin{eqnarray*}
 V_{m}^2&=& \sum_{J\in m}\sum_{i=1}^{r}\left\vert J\right\vert \overline{s}_{m}\left( J,i\right) \log ^{2}\left[ \frac{\bar{s}_{m}\left( J,i\right) }{\hat{s}_{m}\left( J,i\right) }\right]. \nonumber
\end{eqnarray*}
Finally, the set $\Omf$ defined on a minimal partition $m_f$ is
defined in this context by
\begin{eqnarray*}
 \Omf &=&\bigcap\limits_{J\in m}\bigcap\limits_{i=1}^{r}\left\{ \left\vert N_{J}\left( i\right) -\sum_{t\in J}s\left( t,i\right) \right\vert
 \leq \varepsilon \sum_{t\in J}s\left( t,i\right) \right\}.
\end{eqnarray*}

All those quantities are defined in the same manner than in our
general approach, as we can recover $E^i_J=\sum_{t\in J}s\left(
t,i\right)$ and $T^i_J = N_{J} (i)$. The main difference is that
the dimension differs since $r=d+1$ terms are considered in each
sum.

\paragraph{\textbf{Control of $K\left(\bar{s}_{m},\hat{s}_{m}\right)$ by $V_{m}^{2}$}}

\begin{lmm}
 \label{ineg_K3} For all positive densities $p$ and $q$ with respect to $\mu $, one has
 \begin{equation*}
   \frac{1}{2}\int f^{2}\text{ }\left( 1\wedge e^{f}\right) \text{ }p\text{ } d\mu \leq K\left( p,q\right) \leq \frac{1}{2}\int f^{2}\text{
   }\left( 1\vee e^{f}\right) \text{ }p\text{ }d\mu
 \end{equation*}
 if one notes $f=\log \left( \frac{q}{p}\right) $.
\end{lmm}
Applying this lemma (for which a proof can be found  in
\cite{castellan1999modifiedb}), we have on $\Omf$
\begin{eqnarray*}
 K\left( \bar{s}_{m},\hat{s}_{m}\right)  &\geq &\frac{1}{2}\sum_{J\in m}\sum_{i=1}^{r}\left\vert J\right\vert \overline{s}_{m}\left( J,i\right) \left[ 1\wedge \frac{\hat{s}_{m}\left( J,i\right)}{\bar{s}_{m}\left( J,i\right) }\right] \log ^{2}\left[ \frac{\bar{s}_{m}\left( J,i\right) }{\hat{s}_{m}\left( J,i\right) }\right]  \\
 &\geq &\frac{1-\varepsilon \text{ }}{2}V_{m}^2,
\end{eqnarray*}%
and%
\begin{eqnarray*}
 K\left( \bar{s}_{m},\hat{s}_{m}\right)  &\leq &\frac{1}{2}\sum_{J\in m}\sum_{i=1}^{r}\left\vert J\right\vert \overline{s}_{m}\left( J,i\right) \left[ 1\vee \frac{\hat{s}_{m}\left( J,i\right) }{\bar{s}_{m}\left( J,i\right) }\right] \log ^{2}\left[ \frac{\bar{s}_{m}\left( J,i\right) }{ \hat{s}_{m}\left( J,i\right) }\right]  \\
 &\leq &\frac{1+\varepsilon }{2}\text{ }V_{m}^2,
\end{eqnarray*}
resulting in
\begin{equation} \label{RelationVl}
 \frac{1-\varepsilon \text{ }}{2} V_{m}^{2} \leq K\left( \bar{s}_{m},\hat{s}_{m}\right) \leq  \frac{1+\varepsilon\ }{2} V_{m}^{2}\text{,}
\end{equation}

\paragraph{\textbf{Exponential bound for $\chi _{m}^{2}$}}

We want to control the chi-square statistic around its expectation
by using Bernstein's inequality as in Section~\ref{exp}. Here the
required control of $N_{J}(i)$ can be obtained through a direct
application of the bounded version of Bernstein since
$\mathbf{1}_{\{Y=i\}}$ is bounded by $1$. We get:
\begin{eqnarray*} 
 P\left[ \left|\sum_{t\in J}\overline{\mathbf{1}_{\{Y=i\}}}\right|\geq x\right]  &\leq &2\exp \left( -\frac{x^{2}}{2\left( \frac{x}{3}+\sum_{t\in J}s\left( t,i\right) \right) }\right),
\end{eqnarray*}
which, as in the general approach, can be used to obtain
\begin{equation*}
 P\left( \Omf^{c}\right) \leq\frac{C\left( \Gamma ,\rho ,a,r,\varepsilon \right) }{n^{a}}\text{.}
\end{equation*}

Then, since
\begin{equation*}
 \Eb\left[ \chi^2(m,i)\right] =\sum_{J\in m} \frac{Var \left[ \sum_{t\in J} \mathbf{1}_{\{Y=i\}} \right] }{\sum_{t\in J}s\left( t,i\right) }  =\sum_{J\in m}\frac{\sum_{t\in J}s\left( t,i\right) \left[ 1-s\left( t,i\right) \right] }{\sum_{t\in J}s\left( t,i\right) }.
\end{equation*}
and $\forall t,i$ $\rho \leq s\left( t,i\right) \leq 1$, and $r-1
\geq r/2$ for $r \geq 2$, we have the following bounds
\begin{equation}
 \rho \frac{r}{2}|m|\leq \sum_{i=1}^{r}\mathbb{E}_{s}\left[ \chi^2(m,i) \right] \leq r|m|\text{.}
 \label{eq:esperance_X3}
\end{equation}

By controlling the moment of $Z_J (i)$ (as in Section~\ref{exp}) and
using Bernstein's inequality again, we conclude to
\begin{eqnarray*}
 P\left( \sum_{i=1}^{r}\chi^2(m,i) \mathbf{1}_{\Omf }\text{ }\geq r|m|+8r\left( 1+\frac{\varepsilon }{3}\right) \sqrt{x|m|}+4r\left( 1+\frac{\varepsilon}{3}\right)x\right) &\leq &r\exp \left( -x\right).
\end{eqnarray*}

\paragraph{\textbf{Control of $K\left(\bar{s}_{m},\hat{s}_{m}\right)$ by  $\chi _{m}^{2}$}}

The term $V^2_m$ can be written as follows:
\begin{equation*}
 V^2_m= \sum_{J\in m}\sum_{i=1}^{r}|J|\frac{\left[ \hat{s}_{m}(J,i)-\bar{s}_{m}(J,i)\right] ^{2} }{\bar{s}_{m}(J,i)}\left[ \frac{\log \left[ \hat{s}_{m}(J,i)/\bar{s}_{m}(J,i) \right] }{\hat{s}_{m}(J,i)/\bar{s}_{m}(J,i)-1}\right] ^{2},
\end{equation*}
so that using $\frac{1}{1\vee x}\leq \frac{\log {x}}{x-1}\leq
\frac{1}{1\wedge x}\text{ for all }x>0$, on the set $\Omf$, we have:
\begin{equation*}
 \frac{1}{\left( 1+\varepsilon \right) ^{2}}\sum_{i=1}^{r}\chi^2(m,i) \leq V_{m}^{2}  \leq \frac{1}{\left( 1-\varepsilon \right) ^{2}}\sum_{i=1}^{r}\chi^2(m,i).
\end{equation*}
Combining this equation with relation (\ref{RelationVl}) gives, on
$\Omf$
\begin{equation}  \label{chi_l3}
 \frac{1-\varepsilon }{2\left( 1+\varepsilon \right)^{2}}\sum_{i=1}^{r}\chi^2(m,i) \leq K\left(\bar{s}_{m},\hat{s}_{m}\right) \leq \frac{1+\varepsilon }{2\left(1-\varepsilon \right) ^{2}}\sum_{i=1}^{r}\chi^2(m,i).
\end{equation}

\subsubsection{Control of the two other terms}
We explicit here the control of $\bar{\gamma}\left(
\bar{s}_{m}\right) -\bar{\gamma}\left( s\right) $ and
$\barga(\hatsmp)-\barga(\barsmp)$ which appeared in the negligible
constant of the risk.

\begin{itemize}
 \item Control of the term $\bar{\gamma}\left( \bar{s}_{m}\right) -\bar{\gamma}\left( s\right) $. We have to control
   \begin{equation*}
     \mathbb{E}_{s}\left[ \left( \bar{\gamma}\left( \bar{s}_{m}\right) -\bar{\gamma}\left( s\right) \right) \mathbf{1}_{\Omf\left( \varepsilon
\right) }\right] =-\mathbb{E}_{s}\left[ \left( \bar{\gamma}\left(
\bar{s}_{m}\right) -\bar{\gamma}\left( s\right) \right)
\mathbf{1}_{\Omega _{m}\left( \varepsilon \right) ^{c}}\right]
\text{,}
   \end{equation*}%
   and we bound this expectation by
   \begin{eqnarray*}
     \left\vert \Eb\left[ \left( \bar{\gamma}\left(\bar{s}_{m}\right) -\bar{\gamma}\left( s\right) \right) \mathbf{1}_{\Omf}\right] \right\vert &\leq &\left\vert \mathbb{E}\left[ \left( \bar{\gamma}\left( \bar{s}_{m}\right) -\bar{\gamma}\left( s\right) \right) \mathbf{1}_{\Omf^{c}\left( \varepsilon \right) }\right]\right\vert \\
     &\leq &nr\log \left( \frac{1}{\rho }\right) P\left( \Omf^{c}\right)\\
     &\leq & \frac{C\left( \Gamma ,\rho ,a,r,\varepsilon \right)}{n^{a-1}}
   \end{eqnarray*}%

 \item Control of term $\barga(\hatsmp)-\barga(\barsmp)$. We use the following proposition
   \begin{prpstn} \label{K_H3}
     For every $u\in \mathcal{S}$ and any positive $x$,
     \begin{equation*}
       P\left( \bar{\gamma}\left( s\right) -\bar{\gamma}\left( u\right) \geq K\left( s,u\right) -2h^{2}\left( s,u\right) +2xr\right) \leq r\exp \left( -x\right) ,
     \end{equation*}
     where
     \begin{equation*}
       \bar{\gamma}\left( s\right) -\bar{\gamma}\left( u\right) =\sum_{i=1}^{r}\sum_{t=1}^{n}\overline{\mathbf{1}_{\left\{ Y_{t}=i\right\} }}\log \left[ \frac{u\left( t,i\right) }{s\left( t,i\right) }\right] ,
     \end{equation*}
     and $h^{2}$ is the Hellinger distance.
   \end{prpstn}

The proof relies on the same arguments than for  proposition
\ref{s-u}.
\end{itemize}

\subsubsection{Proof of theorem \ref{theo_categorial_direct}}

The proof of this theorem is obtained by following the lines of
Section \ref{prooftheo}: we introduce sets of large probability,
$\OO1$ and $\mO2$ and gather the previous results on the control of
each terms of the main decomposition. This results in the main
inequality :

\begin{equation*}
\Eb \left[ h^{2}\left( s,\hatsmh\right) \right] \leq
\frac{2\text{{}}(2\lambda_{\varepsilon})^{1/3}}{(2\lambda_{\varepsilon})^{1/3}\
-1}\left\{ K\left( s,\barsm\right) +{pen}\left(m\right) \right\}
+C\left(\Sigma ,r,\lambda_{\varepsilon},\rho ,\Gamma \right).
\end{equation*}%

\subsubsection{Oracle-type inequality}

The following proposition gives the risk of an estimator $\hat{s}_m$
for any $m \in \mathcal{M}_n$.
\begin{prpstn}\label{prisk3}
 Under the assumptions of theorem \ref{theo_categorial_direct}, we have
 \begin{eqnarray} \label{eq:Risk-control-direct}
   \Eb\left[ K\left( s,\widehat{s}_{m}\right)\right] &\geq &K(s,\bar{s}_{m})+\frac{1-\varepsilon }{4\left( 1+\varepsilon \right) ^{2}}\ \rho \text{ }r\text{ }|m| -\frac{C\left( \Gamma ,\rho ,a,r,\varepsilon \right) }{n^{a-1}},
 \end{eqnarray}%
 where $a>2$ and $C\left( \Gamma ,\rho ,a,r,\varepsilon \right) $ is a positive constant only depending on $\Gamma $, $a$, $\rho $, $r$ and $\varepsilon $.
\end{prpstn}

Proof: According to inequalities (\ref{eq:esperance_X3}) and
(\ref{chi_l3}), we get the following lower bound of the risk:
\begin{eqnarray*}
\Eb\left[ K(s,\hat{s}_{m})\right] &=&K\left( s,\bar{s}%
_{m}\right) +\Eb \left[ K(\bar{s}_{m},\hat{s}_{m}) \mathbf{1}_{\Omf}\right] +\Eb \left[ K(%
\bar{s}_{m},\hat{s}_{m}) \mathbf{1}_{\Omf^{c}}\right] \\
&\geq & K \left( s,\bar{s}%
_{m}\right) +\frac{1-\varepsilon }{4\left( 1+\varepsilon \right)
^{2}} \rho r|m|-\frac{C\left( \Gamma ,\rho ,a,r,\varepsilon \right)
}{n^{a-1}}.
\end{eqnarray*}%

Using the results of proposition \ref {prisk3} and theorem
\ref{theo_categorial_direct}, we obtain the following oracle-type
inequality:
\begin{crllr}  \label{CorLikelihood}
 Let $pen:\mathcal{M}_{n}\to\mathbb{R}_+$ be defined by (\ref{pen_generale_likelihood}), and assume $s \geq \rho$ for some $\rho> 0$. Then there exists some constant $C$ such that for an exhaustive search,
 \begin{equation} \label{oracle-direct}
   \Eb \left[ h^{2}\left( s,\hatsmh\right) \right ] \leq C\log \left( n\right) \inf_{m\in \mathcal{M}_{n}}\left\{ K\left( s,%
   \hat{s}_{m}\right) \right\} +C(\Gamma ,\rho ,a,r,\lambda_{\varepsilon},\Sigma) \text{.}
 \end{equation}
\end{crllr}
Remark: $K\left( s,\barsm\right) $ is close to $4h^{2}\left(
s,\barsm\right) $ when $||\log (s/\barsm)||_{\infty }<\infty $ and
under this assumption the above upper bound can be stated as follows
\begin{equation*}
 \Eb \left[ h^{2}\left( s,\hatsmh\right) \right] \leq C\log\left( n\right)  \text{ }\inf_{m\in \mathcal{M}_{n}}%
 \Eb\left[ h^{2}\left( s,\hatsm\right) \right]+C(\Gamma ,\rho ,a,r,\lambda_{\varepsilon},\Sigma) .
\end{equation*}

\subsection{Comparison of the constants}

{First, we aim at comparing the bounds of the risk of
our estimator between the general and direct approach (see equations
(\ref{Risk}) and (\ref{risk-direct})).} In the general approach, the
constant is expressed as
$C_{\beta_\varepsilon}=\frac{1+C_2^{-1}(\varepsilon)}{C_2^{-1}(\varepsilon)}$
where
$C_2(\varepsilon)=\beta_{\varepsilon}=\dfrac{E^{max}\kappa}{2m_\varepsilon}(1+\varepsilon)^3$.
 Using the results
of Section \ref{cat:general}, we have $E^{max}< 1$, $\kappa=2$ and
$m_{\varepsilon}=(E^{min}-\varepsilon)^{d+1}$, so that
$C_2(\varepsilon)$ can be bounded by
$$C_2(\varepsilon) \leq \dfrac{(1+\varepsilon)^3}{(E^{min}-\varepsilon)^r}\leq \left(\dfrac{1+\varepsilon}{E^{min}-\varepsilon}\right)^\alpha\quad \text{where} \ \alpha=\max\{r,3\}.$$
We therefore obtain:
\begin{eqnarray*}
C_{\beta_\varepsilon} \leq
\dfrac{(1+E^{min})\beta_{\varepsilon}^{1/\alpha}}{E^{min}\beta_{\varepsilon}^{1/\alpha}-1}
\end{eqnarray*}
which can further be bounded by
$C_{\beta_\varepsilon}=\dfrac{2\beta_{\varepsilon}^{1/\alpha}}{\frac{\beta_{\varepsilon}^{1/\alpha}}{a+a^2d}-1}$
with $\beta_{\varepsilon}> (a+a^2d)^{\alpha}$. \\
In the direct
approach, the constant is $C_{\lambda_{\varepsilon}
}=\frac{2\text{{}}(2\lambda_{\varepsilon})
^{1/3}}{(2\lambda_{\varepsilon} ) ^{1/3}\text{{}}-1}$ with
$\lambda_{\varepsilon}=\frac{1}{2}\frac{(1+\varepsilon)^3}{(1-\varepsilon)^3}>\frac{1}{2}$.\\

In both cases, $\beta_{\varepsilon}$ and $\lambda_{\varepsilon}$ are the constants from the penalty function that are tuned from the data. We can notice that the general shape of $C_{\beta_{\varepsilon} }$ and $C_{\lambda_{\varepsilon} }$ are the same, except that the power of the penalty constant is directly related to the number of categories in the general case while it is fixed to $3$ in the direct approach.\\
This is compensated by, on one hand, a greater constraint on
${\beta_{\varepsilon} }$ since $a>1$ and thus ${\beta_{\varepsilon}
}>r^r$ while ${\lambda_{\varepsilon} }>1/2$, and on the other hand
by a different multiplicative factor in the penalty function, since
we obtain $pen(m)>\beta_{\varepsilon}d|m|(1+4\sqrt{L_m})^2 $ in the
general approach and
$pen(m)>\lambda_{\varepsilon}(d+1)|m|(1+4\sqrt{L_m})^2 $ in the
direct approach. \\

{Then comparing the oracle inequalities given by
(\ref{oracle}) and (\ref{oracle-direct}) in the general and direct
cases respectively results in comparing the constants $C$ that is:
\begin{itemize}
\item $C_{\beta_\varepsilon} \frac{\beta_\varepsilon}{C(\varepsilon)}$ where $C(\varepsilon)$ (in equation (\ref{eq:Risk-control})) is
$C(\varepsilon)=\frac{1}{2}\zeta E^{min}\dfrac{m_{\varepsilon}}{M_\varepsilon}$ in the general approach, \textit{i.e.} $C(\varepsilon)=\frac{1}{2}E^{min}(1-E^{max}){m_{\varepsilon}}$ 
\item $C_{\lambda_\varepsilon} \frac{\lambda_\varepsilon}{C(\varepsilon)}$ where $C(\varepsilon)=\frac{1-\varepsilon }{4\left( 1+\varepsilon \right) ^{2}} \rho $ (cf equation (\ref{eq:Risk-control-direct})) in the direct one.
\end{itemize}

$C_{\beta_\varepsilon} \beta_\varepsilon$ and $\lambda_\varepsilon C_{\lambda_\varepsilon}$ behaving as previously, one looses at least a constant $\frac{1}{2\zeta (E^{min})^{d+1}}$ between the general and the direct approaches.

\section{Simulation study and application to DNA sequences}\label{appli}

In this section, we apply our estimator in two scenarios. The first
is a simulation study where the observations are sampled from the
exponential distribution with piece-wise constant rate parameter. In
this case the distribution is continuous and the sufficient
statistic is unidimensional (and is the variable itself). The second
is an application to a real data-set on the analysis of DNA
sequences in terms of base-composition. The objective is to apply a
segmentation model in order to find homogeneous regions which can be
related to structural and functional biological regions.  In this
case we model the observations with a categorical distribution with
piece-wise constant proportion parameters. The distribution is
therefore discrete and the sufficient statistic is multivariate
(equal to the indicator of the variable taking a category value).

\subsection{Simulation study with exponential distribution}

We simulated $400$ datasets of length $n=10^5$ for which the number
of segments $K$ was drawn from a Poisson distribution with mean
$\bar{K}=50$, and the $K-1$ change-points were sampled uniformly on
$\{2,\dots,n-1\}$ subject to the constraint that segments had to be
of length at least $10$. We considered $4$ sets of values for the
rate parameters. In all scenarios, odd segments had a rate of $0.01$
while the rate on even segments was chosen randomly with probability
$0.4$, $0.2, 0.3$ and $0.1$ among the values $(0.05,0.1,0.02,0.005)$
for datasets $1$ through $100$, $(0.02,0.05,0.015,0.005)$ for
datasets $101$ through $200$, $(0.015,0.02,0.0125,0.007)$ for
datasets $201$ through $300$, and $(0.015,0.02,0.011,0.008)$ for
datasets $301$ through $400$. This resulted in datasets similar to
those shown in Figure \ref{fig:seg}. \newline

In this framework, the segmentation was performed using the pruned
dynamic programming algorithm
\cite{rigaill_pruned_2010,cleynen2014segmentor3isback} from which we
obtained the optimal segmentation for every $K$ up to $Kmax=200$.
The number of segments was then obtained using the penalty function
proposed in (\ref{penalty_exhaustive}) where the constant
is calibrated using the slope heuristic (see \cite{arlot2009data}).\\

To assess the quality of our results, we introduce, as in
\cite{harchaoui2010multiple}, the quantities
$\mathcal{E}(m_{\hat{s}}||m_s)$ and $\mathcal{E}(m_s||m_{\hat{s}})$
between the partitions associated with the true and estimated
distributions $s$ and $\hat{s}$, where
$$\mathcal{E}(m_A||m_B)=\sup_{b\in m_B}\inf_{a\in m_A}|a-b|. $$
The first quantity assesses how our estimated segmentation is able to recover the true change-points. Intuitively, the segmentation with the largest number of segments will have the  greatest chance of yielding a small value of $\mathcal{E}(m_{\hat{s}}||m_s)$. On the contrary, the second quantity judges how relevant the proposed change-points are compared to the true partition: a segmentation with too many segments will necessarily have change-points far from the true ones. Note that the Hausdorff distance can then be recovered as $\sup\{\mathcal{E}(m_{\hat{s}}||m_s),\mathcal{E}(m_s||m_{\hat{s}}) \}.$ \\

The performance of our procedure is finally assessed \textit{via}:
\begin{itemize}
\item The difference between the true number of segments and the estimated, $\Delta=K-\hat{K}$;
\item $\mathcal{E}(m_{\hat{s}}||m_s)$ and  $\mathcal{E}(m_s||m_{\hat{s}})$,  between our estimator and the true segmentation;
\item The Kullback distance between true and estimated distributions, namely $\sum_t \left(\frac{\hat{\lambda_t}}{{\lambda}_t} -\log \frac{\hat{\lambda_t}}{{\lambda}_t} -1\right)$; and
\item The Hellinger distance between true and estimated distributions, namely $\sum_t \left(1-2\frac{\sqrt{\lambda_t \hat{\lambda}_t}}{\lambda_t+\hat{\lambda}_t}\right)$.
\end{itemize}

To allow a fair assessment of our estimator, the
Haussdorff, Kullback and Hellinger criteria are also computed for
the optimal segmentation for the true number of segments, which we
will denote $\hat{m}_{K}$, and the resulting estimated distribution
$\hat{s}_{\hat{m}_{K}}$.
\newline

Our method leads to a tendency to under-estimate the number of
segments, as is shown in Figure \ref{fig:res} a). As in classical
studies of model selection for segmentation, our estimator tends to under-estimate the number of segments in particular when the scenario becomes  more difficult to segment.
Indeed, the number of segments is reduced in order to avoid false
detection. This phenomenon is illustrated in Figure \ref{fig:res} b)
as our estimator (represented through the blue boxplots) yields high
values of $\mathcal{E}(m_{\hat{s}}||m_s)$ due to the missed segments
(subfigure $b_1$) but low values of $\mathcal{E}(m_s||m_{\hat{s}})$
as the segments we propose tend to correspond to true segments
(subfigure $b_2$).
 \begin{figure}[h!]
   \centering
   \includegraphics[width=\textwidth]{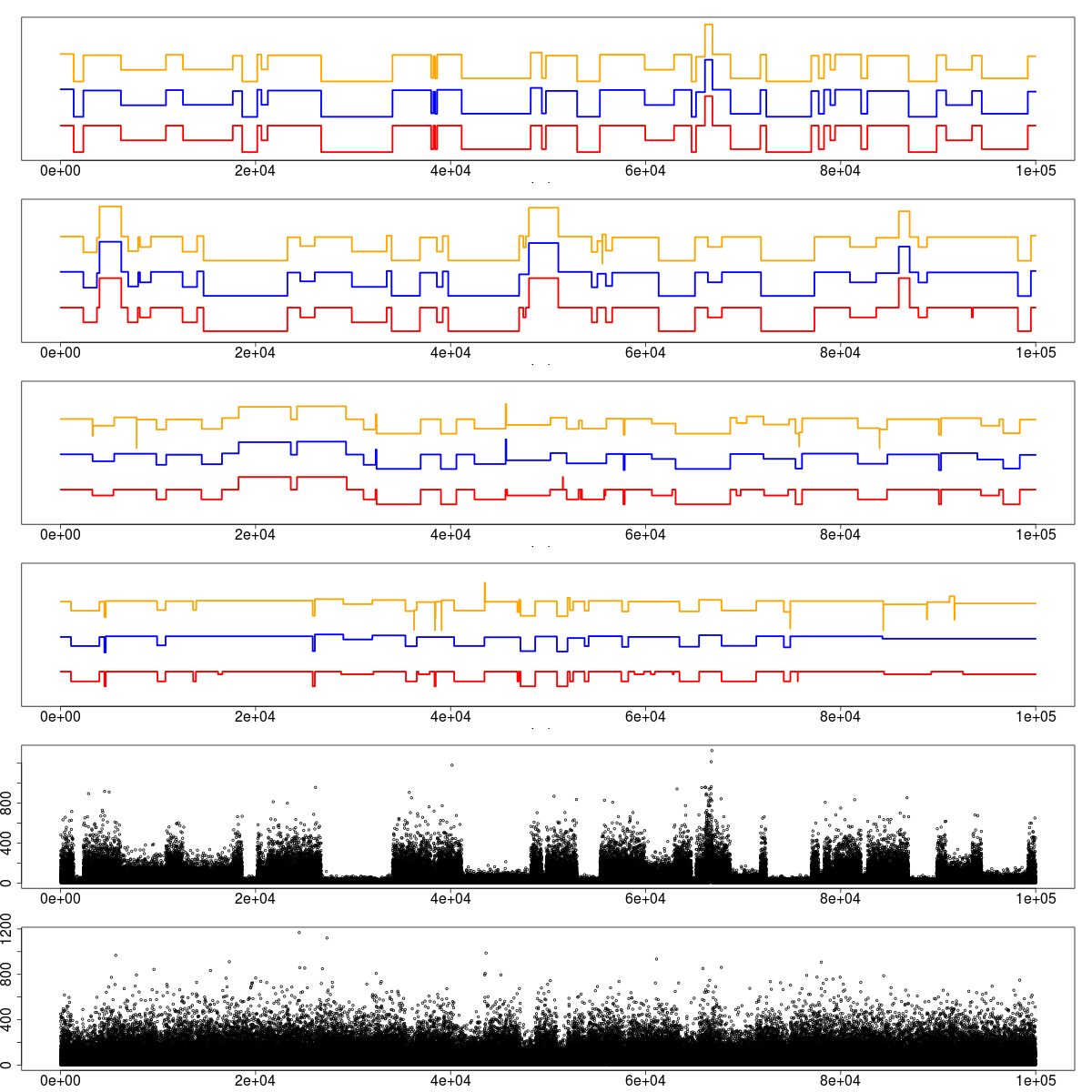}
   \caption{\textbf{Simulated datasets and segmentation.} The four first figures show examples of simulation design for each of the four groups of rate values; the last figure shows an example of dataset to segment in the last group of values. Red lines (bottom) indicate the true distribution (the inverse of the rate, $1/\lambda$, is represented) while blue lines (middle) indicate the estimated distribution and the orange lines (top) indicate the optimal segmentation (w.r.t. likelihood loss) for the true number of segments.}\label{fig:seg}
 \end{figure}
On the contrary, the segmentation $\hat{s}_{\hat{m}_K}$ corresponding to the true number of segments (yellow boxplots) has lower values of $\mathcal{E}(m_{\hat{s}}||m_s)$ as it has more change-points thus lower distances to the missed ones, but higher values of $\mathcal{E}(m_s||m_{\hat{s}})$ as some of its segments are spurious. On a particular example, presented in Figure \ref{fig:seg} by comparing the red line (true segmentation) to the blue one (corresponding to our estimator), and more so on the bottom subfigure illustrating a simulation from the fourth group, we observe that our method tends to fail to recover small segments with rate very close to surrounding ones, whereas the optimal segmentation (as represented by the yellow line) still fails to recover those small segments, thus proposing additional ones, typically very short (average length=9) with very different rate values. \\
Finally, Figures \ref{fig:res} c) and d) show that in terms of
Kullback-Leibler and Hellinger distances, our estimator performs at
least as well than the one with the true number of segments, and
significantly better when the profiles become more difficult to
segment.

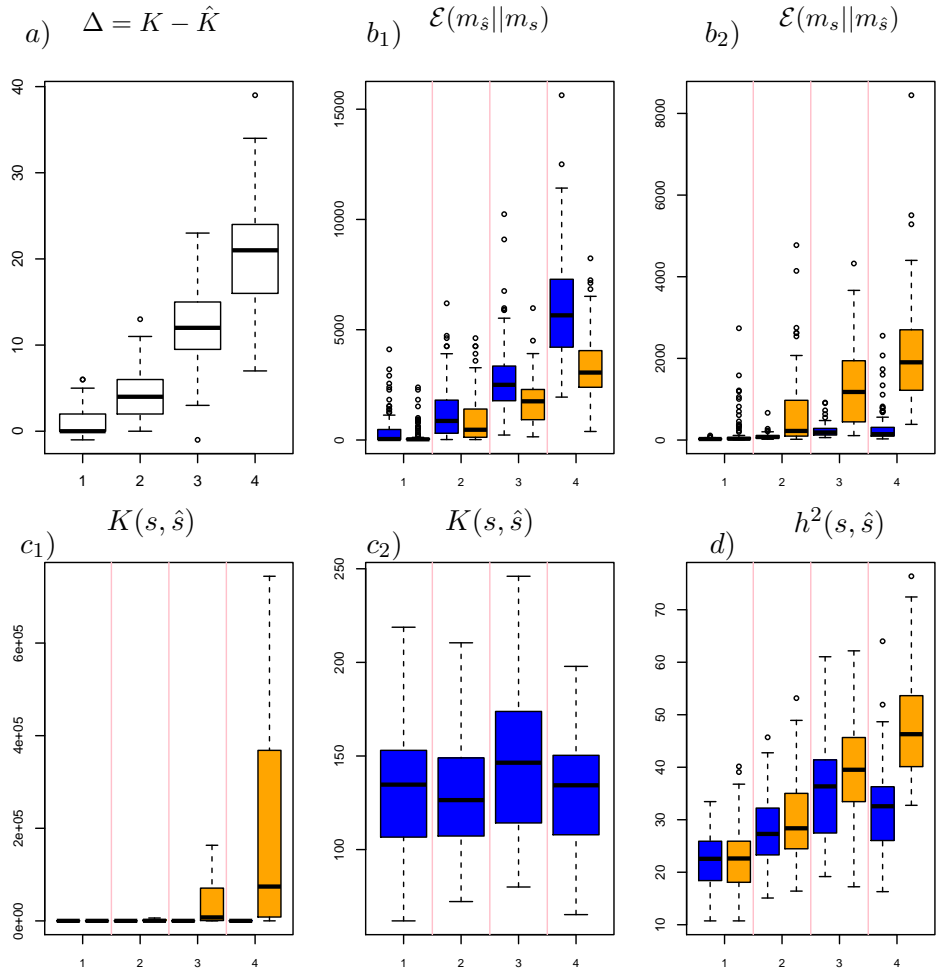
\begin{figure}[h!]
 \centering
\begin{tikzpicture}
  \node at (0,0) {\includegraphics[width=0.8\textwidth]{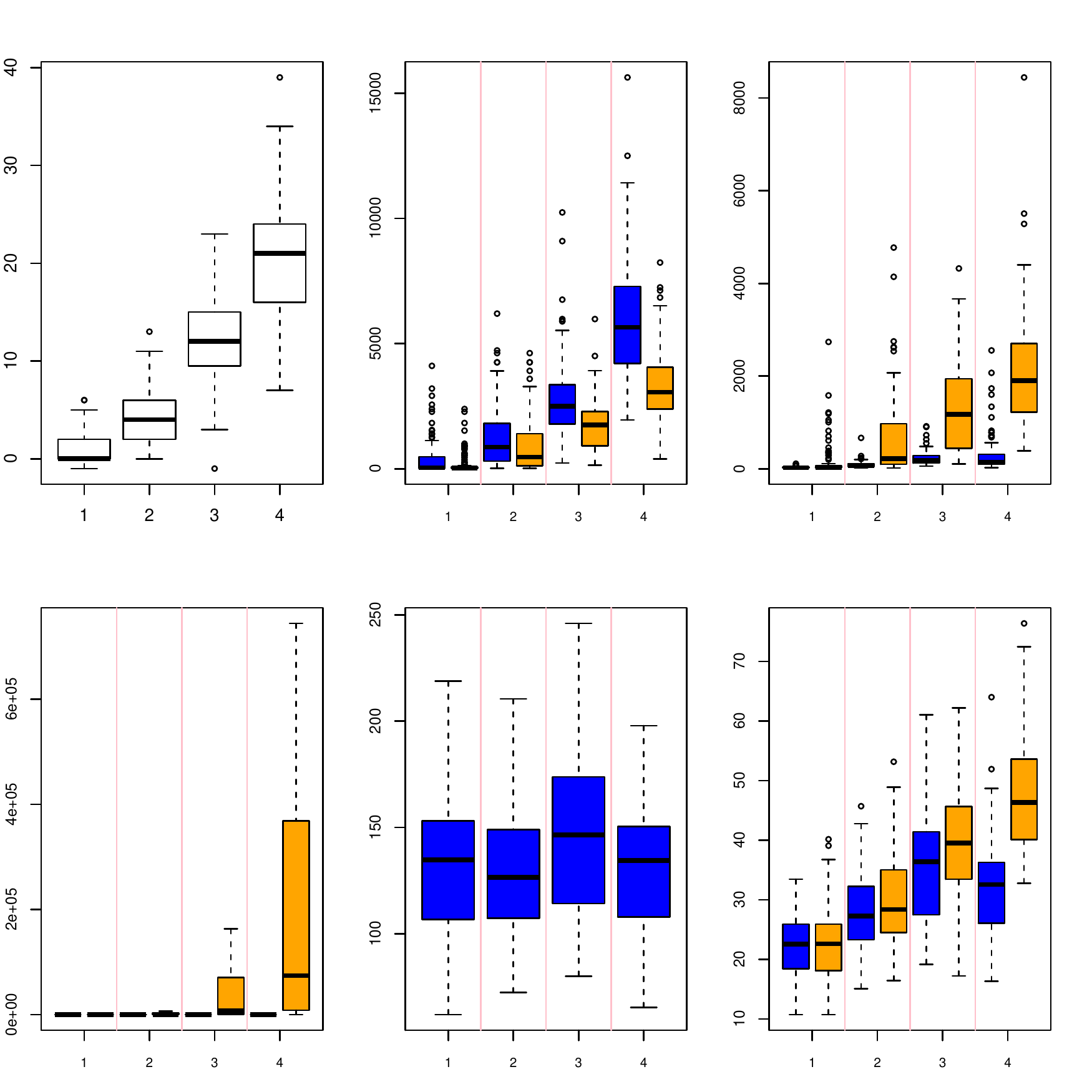}};
  \node at (-6,6.3) {\large{$a)$}};
  \node at (-1.4,6.3) {\large{$b_1)$}};
  \node at (3.1,6.3) {\large{$b_2)$}};
  \node at (-6,-0.5) {\large{$c_1)$}};
  \node at (-1.4,-0.5) {\large{$c_2)$}};
  \node at (3.1,-0.5) {\large{$d)$}};
  \node at (-4.5,6.5) {{$\Delta=K-\hat{K}$}};
  \node at (0,6.5) {{$\mathcal{E}(m_{\hat{s}}||m_s)$}};
  \node at (4.65,6.5) {{$\mathcal{E}(m_s||m_{\hat{s}})$}};
  \node at (-4.5,-0.2) {\large{$K(s,\hat{s})$}};
  \node at (0,-0.2) {\large{$K(s,\hat{s})$}};
  \node at (4.65,-0.2) {\large{$h^2(s,\hat{s})$}};
\end{tikzpicture}
\caption{\textbf{Evaluation of the performance criteria in the four
simulation study frameworks.} a) difference between the true number
of segments and the estimated. Our method tends to underestimate the
number of segments, with performance degradating as the segmentation
difficulty increases. $b_1)$ boxplot of the values of
$\mathcal{E}(m_{\hat{s}}||m_s)$  and $b_2)$
$\mathcal{E}(m_s||m_{\hat{s}})$ over the hundred simulations in each
framework. In each case, left boxplots (blue) assess our estimator,
right boxplots (orange) assess the estimator corresponding to the
true number of segments. $c_1)$ and $c_2)$ Boxplot of the
Kullback-Leibler distance to the true distribution for the estimated
distribution (blue), and the optimal distribution in $K$ segments
(orange) in each simulation framework. In $c_2)$, only
$\hat{s}_{\hat{K}}$ is shown. $ d)$ Same as $c_1)$ but for the
Hellinger distance.} \label{fig:res}
\end{figure}

\subsection{Application of a DNA sequence} \label{sec:Appli}

The objective of this application is to find regions of a DNA
sequence which are homogeneous in terms of  base composition, that
is which present a stability in the frequencies of the four
nucleotide letters. These regions are thought to  correspond to
areas of the genome which are biologically significant. To this end,
we apply our procedure modeling the data with categorical variables
with $d=3$ (see Section \ref{categorical} for the model). \\
Here we consider the bacteriophage Lambda genome with length
$n=48502$ base pairs which is a parasite of the intestinal bacterium
{\it Escherichia coli}. This genome has been used for the comparison
of segmentation methods (see \cite{braun1998statistical} and
references therein) such as HMM (\cite{BA2004}, \cite{Muri98}) or
penalized quasi-likelihood \cite{braun2000multiple}. The data and
its annotation are publicly available from the National Center for
Biotechnology Information (NCBI) pages
at the url adress http://www.ncbi.nlm.nih.gov/. \\
From a computational point of view, the large size of the Lambda
genome hampers the direct use of the classical Dynamic Programming
(DP) algorithm. Here we propose  a hybrid algorithm that consists in
first selecting a small number of relevant change-points using the
CART algorithm \cite{breiman_cart}, and then using dynamic
programming on this set of candidate change-points. As in the
simulation study, the penalty constant is calibrated using the slope
heuristic (see \cite{arlot2009data}). Four change-points
(\textit{i.e.} five segments) are selected by our criterion at
positions $22546$, $27829$, $38004$ and $46528$. The associated
regions are characterized by different base composition as shown in
Figure \ref{fig:lambda_res} which represents the estimated
probabilities of
each base for the obtained segmentation. \\

These change-points are very close (and even on some occasion  the
precise same) as the one obtained in \cite{braun2000multiple}, which
concluded to $3$ more change-points. This reference also supposes
bases to be independent and uses a penalized contrast procedure to
perform the segmentation, and is in this sense the closest approach
to ours. The segments we identify reflect changes in transcription
direction. Indeed, this direction is forward up to base $20855$, it
then switches to reverse from base $22686$ to $37940$, switches back
to forward between $38041$ to $46427$  and finally reverse again
from $46459$ to the end. Note that a refinement has been obtained
when assuming a dependence relationship between bases (see
\cite{BA2004} and \cite{Muri98}).

\begin{figure}[h!]
\centering
\includegraphics[scale=0.8]{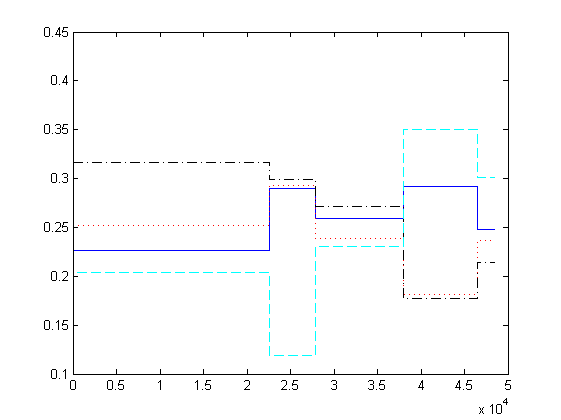}
\caption{Estimated probabilities on each segment of the selected
segmentation: blue and '-' for the base adenosine A, red and ':' for
the base cytosine C, black and '-.' for the base guanine G and cyan
and '--' for the base thymine T.}\label{fig:lambda_res}
\end{figure}

\section{Conclusion}

We have proposed a general approach to the selection of the number
of segments in the general framework where  the data can be modeled
using a distribution from the exponential family. As expected, the
log-partition function and its many properties are instrumental in
the derivation of the bounds and the obtention of
the oracle inequality.\\
The main drawback of our approach is that it can only be applied to
distributions for which the sufficient statistics can be guaranteed
to belong to a positive set. This is mainly due to the use of the
chi-square statistic which was initially defined for the analysis of
count (and thus positive) data. It is very likely that our result
could be extended to overcome this issue by modifying decomposition
\ref{decomposition} and controlling the fluctuation of the
statistics around their expectation using different concentration
inequalities. However, the main goal of this work, beside providing
a general penalty function for model selection, was to underline the
role of the log-partition function in our previous work. Similarly,
this work should easily be extended to multivariate distributions
from
the exponential family. \\
Here, using the particular case of categorical variables as a
example, we have shown that the loss in tightness of the main
constant is not a drastic issue by comparing the results obtained
from the general approach to that from the direct one. \\

Moreover, we have shown in many examples through simulation and
application studies (negative binomial and Poisson distributions in
our previous work, exponential distribution in the simulation study
and categorical distribution in the application to DNA sequences)
that our approach is a powerful method to detect significant changes
in the distribution of the data, which can often be related to
phenomenon of interest. It outperforms  existing criteria (see for
instance \cite{CleynenLebarbier2014}) with a behavior that is
expected  in classical studies of segmentation problems: it tends to
under-estimate the number of segments in order to avoid false
detection. \newline

 \textbf{Acknowledgements} The authors would like to thank
Elodie N\'ed\'elec for her help in the categorical study, and
Mahendra Mariadassou and Hoel Queffelec for helpfull discussions on
exponential families and their geometric properties.

\appendix

\section{Computations for classic laws of the exponential family}\label{app}

Because the constants are strongly related to properties of the cumulants, we distinguish distributions with explicit cumulants from those that can be computed through some recurrence property.\\

We first study laws for which the cumulants are given explicitly.
They include the Poisson, exponential, Gaussian distributions,
etc...
\begin{itemize}
\item In the case of the Poisson distribution with parameter $\lambda$, the log-partition function is analytic on $\mathbb{R}$ and all cumulants are equal to $e^\theta$ so we can use $\zeta=\kappa=1$.\\

\item In the Gaussian case $\mathcal{N}(\mu,\sigma^2)$, we require the mean parameter $\mu$ to be strictly positive in order to define the $\chi^2$ statistic.
In this case, the natural parameters are given by $\thetabf=\left(\dfrac{\mu}{\sigma^2},-\dfrac{1}{2\sigma^2}\right) $, and the sufficient statistic is $T(Y)=(Y,Y^2)$, for which the cumulants are given by explicit formulaes through the log-partition fuction $A(\thetabf)=-\frac{\theta_1^2}{4\theta_2}-\frac{1}{2}\log (\theta_2). $\\
In particular, the cumulants of $Y$ are simply $c_1=-\frac{\theta_1}{2\theta_2}$, $c_2=-\frac{1}{2\theta_2}$ and $c_k=0$ for $k\geq 3$, and we have $Var(Y)=\theta_1^{-1}\mathbb{E}(Y)$ so that $\zeta_1=1/\theta_1^{max}$ and $\kappa_1=1/\theta_1^{min}$ work.\\
For the second order statistic $Y^2$,  for $k\geq 2$ the cumulants are given by
$c_k=\left(-\frac{1}{4}\frac{\theta_1^2}{\theta_2}+\frac{1}{2k}\right)\left(-\frac{1}{\theta_2} \right)^k k!$ on $-|\theta_2|<z<|\theta_2|$ and we have $-\frac{1}{\theta_2}\mathbb{E}(Y^2) \leq Var(Y^2) \leq -\frac{2}{\theta_2}\mathbb{E}(Y^2)$. With the previous notations, we can choose $\tau=-1/\theta_2^{min}$ and $\zeta_2=-1/\theta_2^{min}$ and $\kappa_2=-2/\theta_2^{max}$ work.\\
\\

\item The case of the Pareto distribution also has to be reduced to a fixed scale parameter $x_m$ and a shape parameter $\alpha$ smaller than $-1/\log x_m$.  We have $\zeta=\dfrac{1}{\theta^{max}}\dfrac{1}{{\theta^{max}}\log x_m-1}$. The cumulants are then simply obtained as $c_k=(-1)^k\frac{(k-1)!}{\theta^k} $ for $k>1$ and we can take $\kappa=\frac{\log x_m}{\theta^{min}\log x_m-1}$. \\

\item The Gamma distribution $\Gamma(\alpha,\beta)$ with known shape parameter $\alpha$ has an analytic log-partition function with radius $|\theta|$  where $\theta$ is the {opposite} of the scale parameter $\beta$. The sufficient statistic $Y$ has cumulants given by $c_k=(-1)^k(k-1)!\alpha\theta^{-k}$ and therefore $\tau$ can be chosen as $-\frac{1}{\theta^{max}} $
and finally $\zeta=-1 / \theta^{min} $ and $\kappa=-1/\theta^{max}$. \\ The inverse Gamma distribution with known shape parameter yields the same result since its sufficient statistic is $1/Y$ which, by definition, follows a Gamma distribution.\\

\item The centered Weibull distribution with shape parameter $k$ and scale parameter $\lambda$  belongs to the exponential family provided the shape parameter $k$ is known. The sufficient statistic is then $Y^k$, which, when normalized by $\lambda$, follows an exponential distribution with parameter 1. Results on $\zeta$ and $\kappa$ can therefore easily be deducted from the exponential study. \\

\item The Laplace distribution with known mean positive parameter $\mu$ (positivity is assumed for the definition of $\chi^2$) and scale parameter $b$ has sufficient statistic $|Y-\mu|$ and natural parameter $\theta=-1/b$. As the sufficient statistic follows an exponential distribution with parameter $-\theta$, all results are easily deducted from this study.
\end{itemize}

\vspace{0.5cm}
We then study distributions for which cumulants can be obtained by some recurrence property. This includes the Bernoulli, binomial  (with known parameter $n$), negative binomial  (with known overdispersion parameter  $\phi$), etc... While the recurrence are not derived properly in this section, a complete example is detailed in the application to categorical variables in Section \ref{categorical}.\\
\begin{itemize}
\item The Bernoulli and binomial distributions with fixed number of trials $n$ can be treated together by considering $n=1$ in the first case. The sufficient statistic is $Y$ and the natural parameter is $\theta=\log\frac{p}{1-p}$. The log-partition function, $A(\theta)= n\log\left(1+e^\theta \right)$ is analytic with radius of convergence $\log{(1+\frac{1+e^\theta}{e^\theta})}$ and we can take $R=\log 2$. While the cumulants are not given by an explicit formula, they can be obtained with the recurrence : $c_1=np$ and $c_{k+1}=p(1-p)\frac{\delta c_{k}}{\delta p}$ for $k\geq 1$, which is the same as in the categorical study. Similar computations lead to the results of Table \ref{usual}.\\

\item Similarly, the negative binomial distribution belongs to the exponential family provided the overdispersion parameter $\phi$ is known. The natural parameter is then $\theta=\log{(1-p)}$ and the sufficient statistic is $Y$. Cumulants verify the recurrence property $c_1=\phi \frac{1-p}{p} $ and $c_{k+1}=(1-p) \frac{\delta c_{k}}{\delta (1-p)}$ and we can easily show that $c_{k}=\phi\frac{1-p}{p^{k}}f_k(p)$ with $f_k(p)$ a polynomial in $p$ with degree at most $k-2$ and sum of absolute coefficients less than $2^{k-3}k!$. This leads to the choices $\zeta=1$ and $\kappa=e^{-\theta_{min}}$.
\end{itemize}

\bibliographystyle{bmc_article}
\bibliography{Biblio}
\end{document}